\documentstyle[amstex,amssymb,11pt]{article}

\begin{document}
    \pagestyle{plain} \setlength{\baselineskip}{1.3\baselineskip}
    \setlength{\parindent}{\parindent}

\title{{\bf Coordinates on Schubert cells, Kostant's harmonic forms,
and the Bruhat-Poisson structure on $G/B$}}
\author{Jiang-Hua Lu 
\thanks{Research partially supported by
an NSF Postdoctorial Fellowship.} \\ Department of
 Mathematics, University of Arizona, Tucson, AZ 85721
 USA}
 \maketitle

\newtheorem{thm}{Theorem}[section]
\newtheorem{lem}[thm]{Lemma}
\newtheorem{prop}[thm]{Proposition}
\newtheorem{cor}[thm]{Corollary}
\newtheorem{rem}[thm]{Remark}
\newtheorem{exam}[thm]{Example}
\newtheorem{nota}[thm]{Notation}
\newtheorem{dfn}[thm]{Definition}
\newtheorem{ques}[thm]{Question}
\newtheorem{eq}{thm}

\newcommand{\rw}{\rightarrow}
\newcommand{\lrw}{\longrightarrow}
\newcommand{\rhu}{\rightharpoonup}
\newcommand{\lhu}{\leftharpoonup}
\newcommand{\Map}{\longmapsto}
\newcommand{\qed}{\begin{flushright} {\bf Q.E.D.}\ \ \ \ \
                  \end{flushright} }
\newcommand{\beqa}{\begin{eqnarray*}}
\newcommand{\eeqa}{\end{eqnarray*}}

\newcommand{\la}{\mbox{$\langle$}}
\newcommand{\ra}{\mbox{$\rangle$}}
\newcommand{\lala}{\mbox{$\la \! \la ~,~\ra \!\ra$}}

\newcommand{\ot}{\mbox{$\otimes$}}
\newcommand{\xa}{\mbox{$x_{(1)}$}}
\newcommand{\xb}{\mbox{$x_{(2)}$}}
\newcommand{\xc}{\mbox{$x_{(3)}$}}
\newcommand{\ya}{\mbox{$y_{(1)}$}}
\newcommand{\yb}{\mbox{$y_{(2)}$}}
\newcommand{\yc}{\mbox{$y_{(3)}$}}
\newcommand{\yd}{\mbox{$y_{(4)}$}}
\renewcommand{\aa}{\mbox{$a_{(1)}$}}
\newcommand{\ab}{\mbox{$a_{(2)}$}}
\newcommand{\ac}{\mbox{$a_{(3)}$}}
\newcommand{\ad}{\mbox{$a_{(4)}$}}
\newcommand{\ba}{\mbox{$b_{(1)}$}}
\newcommand{\bt}{\mbox{$b_{(2)}$}}
\newcommand{\bc}{\mbox{$b_{(3)}$}}
\newcommand{\ca}{\mbox{$c_{(1)}$}}
\newcommand{\cb}{\mbox{$c_{(2)}$}}
\newcommand{\cc}{\mbox{$c_{(3)}$}}
\newcommand{\calH}{\mbox{$\cal H$}}
\newcommand{\calS}{\mbox{$\cal S$}}
\newcommand{\hyperH}{\mbox{$\Bbb H$}}
\newcommand{\boldC}{\mbox{$\Bbb C$}}

\newcommand{\ts}{\mbox{$\sigma$}}
\newcommand{\las}{\mbox{${}_{\sigma}\!A$}}
\newcommand{\lasone}{\mbox{${}_{\sigma'}\!A$}}
\newcommand{\ras}{\mbox{$A_{\sigma}$}}
\newcommand{\rds}{\mbox{$\cdot_{\sigma}$}}
\newcommand{\lds}{\mbox{${}_{\sigma}\!\cdot$}}

\newcommand{\bb}{\mbox{$\bar{\beta}$}}
\newcommand{\bg}{\mbox{$\bar{\gamma}$}}

\newcommand{\id}{\mbox{${\em id}$}}
\newcommand{\Fun}{\mbox{${\em Fun}$}}
\newcommand{\End}{\mbox{${\em End}$}}
\newcommand{\Hom}{\mbox{${\em Hom}$}}
\newcommand{\ta}{\mbox{${\mbox{$\scriptscriptstyle A$}}$}}
\newcommand{\ms}{\mbox{${\mbox{$\scriptscriptstyle M$}}$}}
\newcommand{\ap}{\mbox{$A_{\mbox{$\scriptscriptstyle P$}}$}}
\newcommand{\tx}{\mbox{$\mbox{$\scriptscriptstyle X$}$}}
\newcommand{\pp}{\mbox{$\pi_{\mbox{$\scriptscriptstyle P$}}$}}
\newcommand{\pg}{\mbox{$\pi_{\mbox{$\scriptscriptstyle G$}}$}}
\newcommand{\asemi}{\mbox{$\ap \#_{\sigma} A^*$}}
\newcommand{\dsemi}{\mbox{$A \#_{\Delta} A^*$}}

\newcommand{\semi}{\mbox{$\times_{{\frac{1}{2}}}$}}
\newcommand{\fk}{\mbox{${\frak k}$}}
\newcommand{\fa}{\mbox{${\frak a}$}}
\newcommand{\fd}{\mbox{${\frak d}$}}
\newcommand{\ft}{\mbox{${\frak t}$}}
\newcommand{\fg}{\mbox{${\frak g}$}}
\newcommand{\fb}{\mbox{${\frak b}$}}
\newcommand{\fh}{\mbox{${\frak h}$}}
\newcommand{\fn}{\mbox{${\frak n}$}}
\newcommand{\fp}{\mbox{${\frak p}$}}
\newcommand{\fbp}{\mbox{${\frak b}_{+}$}}
\newcommand{\fbm}{\mbox{${\frak b}_{-}$}}
\newcommand{\fnp}{\mbox{${\frak n}_{+}$}}
\newcommand{\fnm}{\mbox{${\frak n}_{-}$}}
\newcommand{\fgs}{\mbox{${\frak g}^*$}}
\newcommand{\wg}{\mbox{$\wedge {\frak g}$}}
\newcommand{\wgs}{\mbox{$\wedge {\frak g}^*$}}
\newcommand{\wxl}{\mbox{$x_1 \wedge x_2 \wedge \cdots \wedge x_l$}}
\newcommand{\wxk}{\mbox{$x_1 \wedge x_2 \wedge \cdots \wedge x_k$}}
\newcommand{\wyl}{\mbox{$y_1 \wedge y_2 \wedge \cdots \wedge y_l$}}
\newcommand{\wxkm}{\mbox{$x_1 \wedge x_2 \wedge \cdots \wedge x_{k-1}$}}
\newcommand{\wxik}{\mbox{$\xi_1 \wedge \xi_2 \wedge \cdots \wedge \xi_k$}}
\newcommand{\wxikm}{\mbox{$\xi_1 \wedge \cdots \wedge \xi_{k-1}$}}
\newcommand{\wetal}{\mbox{$\eta_1 \wedge \eta_2 \wedge \cdots \wedge \eta_l$}}

\newcommand{\winv}{\mbox{$(\wedge \fg_{1}^{\perp})^{\fg_1}$}}
\newcommand{\wetak}{\mbox{$\eta_1 \wedge \cdots \wedge \eta_k$}}
\newcommand{\gonep}{\mbox{$\fg_{1}^{\perp}$}}
\newcommand{\wonep}{\mbox{$\wedge \fg_{1}^{\perp}$}}

\newcommand{\db}{\mbox{$\fd = \fg \bowtie \fgs$}}
\newcommand{\fds}{\mbox{${\scriptscriptstyle {\frak d}}$}}
\newcommand{\fl}{\mbox{${\frak l}$}}

\newcommand{\Gs}{\mbox{$G^*$}}
\newcommand{\pis}{\mbox{$\pi_{\sigma}$}}
\newcommand{\ea}{\mbox{$E_{\alpha}$}}
\newcommand{\eb}{\mbox{$E_{-\alpha}$}}
\newcommand{\Bm}{\mbox{$ {}^B \! M$}}
\newcommand{\kBm}{\mbox{$ {}^B \! M^k$}}
\newcommand{\Bb}{\mbox{$ {}^B \! b$}}
\renewcommand{\epsilon}{\mbox{$\varepsilon$}}

\newcommand{\cfg}{\mbox{$C(\fg \oplus \fgs)$}}
\newcommand{\ps}{\mbox{$\pi^{\#}$}}
\newcommand{\tpi}{\mbox{$\tilde{\pi}$}}
\newcommand{\backl}{\mathbin{\vrule width1.5ex 
height.4pt\vrule height1.5ex}}
 
\newcommand{\bx}{\mbox{${\bar{x}}$}}
\newcommand{\by}{\mbox{${\bar{y}}$}}
\newcommand{\bz}{\mbox{${\bar{z}}$}}
\newcommand{\pgs}{\mbox{${\pi_{\mbox{\tiny G}^{*}}}$}}

\newcommand{\tp}{\mbox{$\varphi$}}
\newcommand{\sn}{\mbox{$s_{\scriptscriptstyle N}$}}
\newcommand{\tn}{\mbox{$t_{\scriptscriptstyle N}$}}
\newcommand{\sm}{\mbox{$s_{\scriptscriptstyle M}$}}
\newcommand{\tm}{\mbox{$t_{\scriptscriptstyle M}$}}
\newcommand{\en}{\mbox{$\epsilon_{\scriptscriptstyle N}$}}
\newcommand{\mem}{\mbox{$\epsilon_{\scriptscriptstyle M}$}}

\newcommand{\lki}{\mbox{$< \!<$}}
\newcommand{\rki}{\mbox{$>\!>$}}
\newcommand{\Xa}{\mbox{$X_{\alpha}$}}
\newcommand{\Ya}{\mbox{$Y_{\alpha}$}}

\newcommand{\qa}{\mbox{$Q_{A}$}}
\newcommand{\pf}{\mbox{$\tilde{\pi}(df)$}}

\newcommand{\cala}{{\cal A}}
\newcommand{\calb}{{\cal B}}
\newcommand{\calc}{{\cal C}}
\newcommand{\cald}{{\cal D}}
\newcommand{\cale}{{\cal E}}
\newcommand{\calf}{{\cal F}}
\newcommand{\calg}{{\cal G}}
\newcommand{\calh}{{\cal H}}
\newcommand{\cali}{{\cal I}}
\newcommand{\calj}{{\cal J}}
\newcommand{\calk}{{\cal K}}
\newcommand{\call}{{\cal L}}
\newcommand{\calm}{{\cal M}}
\newcommand{\caln}{{\cal N}}
\newcommand{\calo}{{\cal O}}
\newcommand{\calp}{{\cal P}}
\newcommand{\calq}{{\cal Q}}
\newcommand{\calr}{{\cal R}}
\newcommand{\calv}{{\cal V}}
\newcommand{\calx}{{\cal X}}
\newcommand{\caly}{{\cal Y}}

\newcommand{\pdp}{\mbox{$\partial_{\pi}$}}
\newcommand{\tdp}{\mbox{$\tilde{\delta}_{\pi}$}}

\newcommand{\ppt}{\mbox{$\partial^{\operatorname{tot}}$}}
\newcommand{\tot}{\mbox{$\operatorname{Tot}$}}
\newcommand{\iep}{\mbox{$i_{\exp_{\wedge} \pi}$}}
\newcommand{\ienp}{\mbox{$i_{\exp_{\wedge} (-\pi)}$}}

\newcommand{\tc}{\mbox{$\tot ~ {\cal C}$}}
\newcommand{\ttc}{\mbox{$\tilde{\cal C}_{\bullet \bullet}$}}
\newcommand{\tcc}{\mbox{$\tot ~ \tilde{\cal C}$}}

\newcommand{\hd}{\mbox{$H\!D$}}
\newcommand{\dw}{\mbox{$\dot{w}$}}
\newcommand{\dc}{\mbox{$\dot{\gamma}$}}
\newcommand{\dgm}{\mbox{$\dot{\gamma}_{\scriptscriptstyle l-1}$}}
\newcommand{\dg}{\mbox{$\dot{\gamma}_{\scriptscriptstyle l}$}}
\renewcommand{\sl}{\mbox{${\scriptscriptstyle l}$}}
\newcommand{\zc}{\mbox{$\{z_1, ~ \bar{z}_1, ~ z_2, ~ \bar{z}_2, ~ ... ~ ,
z_{\sl}, ~ \bar{z}_{\sl}\}$}}
\newcommand{\sdg}{\mbox{$(n_1 \circ \dot{\gamma_1})
(n_2 \circ \dot{\gamma_2}) \cdots
(n_{\sl} \circ \dot{\gamma_{\sl}})$}}
\newcommand{\sdgm}{\mbox{$(n_1 \circ \dot{\gamma_1})
(n_2 \circ \dot{\gamma_2}) \cdots
(n_{\sl -1} \circ \dot{\gamma}_{\sl-1})$}}
\newcommand{\Nhat}{\mbox{$N_{\hat{\gamma}_{\sl}}$}}
\newcommand{\nhat}{\mbox{$n^{\hat{\gamma}_{\sl}}$}}

\tableofcontents

\section{Introduction and Notation}
\label{sec_intro}

This work grew out of my attempts to understand the relations
between results of Kostant \cite{ks:63} on the de Rham cohomology of 
a flag manifold and Poisson geometry.

\bigskip
Let $X = G/B$ be a flag manifold, where $G$ is a complex semi-simple
Lie group and $B$ is a Borel subgroup of $G$. Let $K$ be a 
compact real form of $G$, so $X = G/B \cong K/T$, where
$T = K \cap B$ is a maximal torus of $K$. In \cite{ks:63}, Kostant 
constructs, for each element $w$ in the Weyl group $W$, an explicit 
$K$-invariant closed differential form $s^{w}$ on $X$ with 
$\deg(s^w) = 2 l(w)$, where $l(w)$ is the length of $w$ (see 
Section \ref{sec_kost-thm}).
The cohomology classes of the $s^w$'s form a basis of $H(X, {\Bbb C})$
that, up to constant multiples,
is dual to the basis of the homology of $X$ formed by the 
closures of the Bruhat (or Schubert) cells in $X$. These $K$-invariant forms 
on $X$ are $(d, \partial)$-harmonic, where $d$ is the de Rham 
differential operator and $\partial$ is a degree $-1$ operator 
introduced by Kostant. 

\bigskip
Our work was first motivated by wanting to 
understand the nature of the operator $\partial$ in terms of
Poisson geometry.

\bigskip
The Poisson structure that is relevant to Kostant's theorem 
is the so-called Bruhat-Poisson structure on $X$
\cite{lu-we:poi}. It has its origin in
the theory of quantum groups. It has the special property
that its symplectic leaves are precisely the Bruhat cells
in $X$ (and hence the name). More properties of the Bruhat-Poisson
structure are reviewed in Section \ref{sec_bruhat}. 

\bigskip
The Bruhat-Poisson structure gives rise to the degree $-1$ operator 
\[
\partial_\pi ~ = ~ i_{\pi} d ~ - ~ d i_\pi
\]
on the space of differential forms on $X$ called the 
Koszul-Brylinski operator \cite{ko:crochet} \cite{by:homo}, where 
$\pi$ is the bi-vector field on $X$ defining this Poisson structure,
and $i_\pi$
denotes the contraction operator of differential froms with $\pi$.
It satisfies $\partial_{\pi}^{2} = 0$. Its
homology is called the Poisson homology of $\pi$. Poisson homology
is closely related to cyclic homology of associative algebras, as 
is shown in \cite{by:homo}. It turns out that 
when restricted to $K$-invariant differential
forms on $X$, Kostant's operator $\partial$ and the Koszul-Brylinski
operator $\partial_\pi$ differ by the contraction operator by 
a vector field $\theta_0$ called the modular vector field of
$\pi$ 
(see \cite{we:modular} \cite{bz:outer} \cite{elw:modular}).
More explicitly, it is 
the infinitesimal
generator of the $K$-action on
$X$ in the direction of the element $2iH_{\rho}$ with
$\rho$ being half of the sum of all positive roots. 
This result, together with others on the Poisson (co)homology
of the Bruhat-Poisson structure, can be found in \cite{e-l:poi}.

\bigskip
In this paper, we relate Kostant's harmonic forms on $X$ with
the Bruhat-Poisson structure.

\bigskip
More precisely, the Bruhat-Poisson structure on $X = K/T$ is 
$T$-invariant (but not $K$-invariant). Thus, 
each Bruhat cell $\Sigma_w$ inherits a $T$-invariant symplectic 
structure $\Omega_w$. Use 
$\phi_w: ~ \Sigma_w \rightarrow \ft^*$
to denote the moment map and let 
$\mu_w$ be 
the Liouville volume form on $\Sigma_w$ defined
by $\Omega_w$. Theorem \ref{thm_main} says that
when restricted to the cell $\Sigma_w$, Kostant's form $s^w$
is related to the Liouville volume form $\mu_w$ by
\[
s^w|_{\Sigma_w} ~ = ~ e^{\la \phi_w, ~ 2iH_{\rho} \ra} \mu_w.
\]
Notice that the function $\la \phi_w, ~ 2iH_{\rho} \ra$ 
is the Hamiltonian function for the modular 
vector field $\theta_0$ on $\Sigma_w$. (The vector field $\theta_0$
is not globally Hamiltonian (see Section \ref{sec_kost-thm}), but
it is on each cell). Theorem \ref{thm_main} thus expresses Kostant's
harmonic forms totally in terms of data coming from the Bruhat-Poisson
structure. In particular, it shows that the integral $\int_{\Sigma_w}
s^w$ is of the Duistermaat - Heckman type. Wanting to see this was
another motivation for this work.  

\bigskip
Theorem \ref{thm_main} is proved by writing everything down
in some coordinates 
\[
\{z_1, ~ \bar{z}_1, ~ z_2, ~ \bar{z}_2, ~ ... ~ , z_{\sl}, ~ \bar{z}_{\sl} \}
\]
on each Bruhat cell $\Sigma_w$, where $l = l(w)$. 
These coordinates are motivated by, but are
independent of,  the Bruhat-Poisson structure.
Among the quantities that we write down explicitly in the 
coordinates are (see the later sections for the notation) 

\begin{itemize}
\item
(Theorem \ref{thm_sym-cor}) the symplectic $2$-form $\Omega_w$ on 
$\Sigma_w$ and thus the Liouville volume form $\mu_w$:
\beqa
\Omega_w & = & \sum_{j = 1}^{l}
{\frac{i}{\ll \alpha_j, ~ \alpha_j \gg}}
~ {\frac{1}{1 + |z_j|^2}} ~dz_j \wedge d \bar{z}_j\\
\mu_w & = & \prod_{j=1}^{l} {\frac{i}{\ll \alpha_j, ~ \alpha_j \gg}}
~{\frac{1}{1 + |z_j|^2}}  ~ dz_j \wedge d \bar{z}_j
\eeqa

\item
(Theorem \ref{thm_sym-cor}) the moment map for the $T$-action 
on $(\Sigma_w, ~ \Omega_w)$:
\[
\phi_w: ~~ \Sigma_w \lrw \ft^*: \hspace{.2in} \phi_w ~ = ~
\sum_{j=1}^l (-{\frac{1}{2}} \log(1+|z_j|^2) \check{H}_{\alpha_j})
\]

\item
(Theorem \ref{thm_haar}) the Haar measure $dn$ on the group $N_w
= N \cap w N_{-} w^{-1}$ that
parametrizes $\Sigma_w$:
\[
dn ~ = ~ \left( \prod_{j=1}^{l} 
{\frac{i \ll \rho, ~ \beta_j \gg}{ \pi \ll \beta_j, ~ \beta_j \gg}}
\right) \prod_{j=1}^{l}
 ~(1 + |z_j|^2)^{{\dfrac{2 \ll \rho, ~ \beta_j \gg}
{\ll \beta_j, ~ \beta_j \gg}} ~ - ~1} dz_j \wedge d \bar{z}_j
\]

\item
(Theorem \ref{thm_a}) the $A$-component $a_w(n)$ in the Iwasawa
decomposition of the element $\dw^{-1} n \dw$ for $n \in N_w$
(where $\dw$ is any representative of $w$ in $K$):
\[
a_w(n) ~ = ~ \prod_{j=1}^{l} \exp({\frac{1}{2}} \log(1 + |z_j|^2)
\check{H}_{\beta_j})
\]

\item
(Theorem \ref{thm_main}) Kostant's harmonic form $s^w$ restricted to 
$\Sigma_w$:
\[
s^w|_{\Sigma_w} ~ = ~ \prod_{j=1}^{l}
{\frac{i }{\ll \alpha_j, ~ \alpha_j \gg}} (1+|z_j|^2)^{
{-~ \dfrac{2 \ll \rho, ~ \alpha_j \gg}{\ll \alpha_j, ~ \alpha_j \gg}} 
~ -~ 1}
dz_j \wedge d \bar{z}_j.
\]
\end{itemize}

\bigskip
By comparing these explicit formulas, we immediately get
(Corollaries \ref{cor_aw} and \ref{cor_liouville-haar})
\beqa
\phi_w & = & Ad_{\dot{w}} \log a_w(n) \\
\mu_w & = & \left(\prod_{j=1}^{l} {\frac{\pi}{\ll
\rho, ~ \beta_j \gg}} \right) a_{w}(n)^{-2 \rho} dn.
\eeqa
These relate the moment map $\phi_w$ and the Liouville volume form
$\mu_w$ to the familiar map $a_w$ and the Haar measure 
$dn$ on $N_w$. This is desirable for understanding the 
geometry of the Bruhat-Poisson structure. The relation between
the differential form $s^w|_{\Sigma_w}$ and the Liouville form
$\mu_w$ also follows
immediately from these formulas.

\bigskip
The coordinates $\zc$ are presented first in 
Section \ref{sec_cor}. Here we derive the formulas 
for the Haar measure $dn$ of $N_w$ and for the map 
$a_w: N \rightarrow A$.
As one other application of these coordinates, we
show that
Harish-Chandra's formula for the $c$-function 
(in the case of a complex group) follows easily as 
a product of $1$-dimensional integrals. 
Our calculation here is easier because we have pushed the 
usual induction argument on the length of $w$ into the calculations for
the explicit formulas for $a_w(n)$ and $dn$. The total 
amount of effort is probably the same.

\bigskip
In Section \ref{sec_bruhat}, we review some of the properties of
the Bruhat-Poisson structure on $K/T$. We derive 
the formulas for the  symplectic form 
$\Omega_w$,  thus also the Liouville measure $\mu_w$,
and the moment map for the $T$ action on each 
Bruhat cell in the $\zc$-coordinates. By using the formulas
for $a_w(n)$ and $dn$ given in Section \ref{sec_cor}, we 
arrive at the (coordinate-free) interpretations for both the moment
map $\phi_w$ and the Liouville measure $\mu_w$
in terms of $a_w(n)$ and $dn$ 
as given in Corollaries \ref{cor_aw} and \ref{cor_liouville-haar}.

\bigskip
Kostant's harmonic forms are reviewed in Section \ref{sec_kost-thm}.
Theorem \ref{thm_main} is given as an easy corollary of our formulas
in the $z$-coordinates. 
Applications of Theorem \ref{thm_main} to the calculations of the
Poisson (co)homology of the Bruhat-Poisson structure are given in
\cite{e-l:poi}.

\bigskip
In the Appendix, we discuss how our coordinates $\zc$ are related
to the (complex) Bott-Samelson coordinates.

\bigskip
We now fix the notation.

\bigskip
Let $G$ be a finite dimensional complex semi-simple Lie group
with Lie algebra $\fg$. Let $H$ be a Cartan subgroup and 
$\fh$ its Lie subalgebra. Denote by $R$ the set of all roots
of $\fg$ with respect to $\fh$. Let $R^{+}$ be a choice of 
positive roots. We will also write $\alpha > 0$ for 
$\alpha \in R^+$. Let $\fb = \fb_{+}$ be the Borel subalgebra
spanned by $\fh$ and all the positive root vectors, and let $B$ be
the corresponding Borel subgroup.

\bigskip
Let $W$ be the Weyl group of $G$ relative to $H$.
The Bruhat decomposition
\[
G ~ = ~ \bigcup_{w \in W} B \dot{w} B,
\]
where $\dot{w}$ is a representative of $w$ in the normalizer of $H$ in
$G$, gives rise to the Bruhat decomposition
\[
G/B ~ = ~ \bigcup_{w \in W} B \dot{w} B / B
\]
of the flag manifold $G/B$ into a disjoint union of cells.
These cells, denoted by $\Sigma_w  = B \dot{w} B / B$, are called 
Bruhat or Schubert cells, and 
their closures, $X_w = \overline{\Sigma_w}$, are called
the Schubert varieties.  
 
\bigskip
We choose a compact real form  $\fk$ of $\fg$ as follows:
let $\ll ~ \gg$ be the Killing form of $\fg$. For each 
positive root $\alpha$, denote by $H_{\alpha}$ the image of $\alpha$
under the isomorphism $\fh^* \rightarrow \fh$ via $\ll ~ , ~ \gg$,
i.e., for any $H \in \fh$,
\[
\ll H_{\alpha}, ~ H \gg  ~ = ~ \alpha(H).
\]
Choose root vectors $\ea$ and $\eb$ for $\alpha$ and 
$-\alpha$ respectively such that $\ll \ea, ~ \eb \gg
~ = ~ 1$.
Then $[\ea, \eb ] = H_{\alpha}$. Set
\begin{equation}
\label{eq_xaya}
\Xa ~ = ~ \ea ~ - ~ \eb, \hspace{.3in} \Ya ~ = ~ i(\ea ~ + ~ \eb).
\end{equation}
The real subspace
\[
\fk ~ = ~ \text{span}_{\Bbb R} \{iH_{\alpha}, ~ \Xa, ~ \Ya: ~ 
\alpha > 0\}
\]
is a compact real form of $\fg$. Let $K$ be the corresponding 
compact subgroup of $G$. The intersection $T = K \cap B$ is
a maximal torus of $K$, and its Lie algebra is 
\[
\ft ~ = ~ \text{span}_{\Bbb R} \{iH_{\alpha}: ~ \alpha > 0\}.
\]
Let $\fa ~ = ~ i \ft$, and let $\fn = \fn_{+}$ be the subalgebra of $\fg$
spanned by all the positive root vectors. Let $A$ and $N$ be the 
corresponding subgroups. Then
\[
G ~ = ~ KAN
\]
is the Iwasawa decomposition of $G$ as a real semi-simple Lie group.
It follows that the projection map from $K \subset G$ to $G/B$ 
induces an isomorphism from $K/T$ to $G/B$. From now on, we will identify
$G/B$ with $K/T$ this way.

\bigskip
For $g \in G$, let $g = k_{g} a_{g} n_{g}$ be the Iwasawa
decomposition of $g$. The map
\begin{equation}
\label{eq_G-on-K}
G \times K \lrw K: ~~ (g, ~ k) \Map k_{g \! k}
\end{equation}
defines a left action of $G$ on $K$. We will denote this
action by $(g, k) \rightarrow g \circ k$. 

\bigskip
{\bf Acknowledgement} The author would like to thank Professor
Kostant for explaining to her many of his results, published or
unpublished, and to Sam Evens
for answering many questions. She would also like to thank 
Viktor Ginzburg and Doug Pickrell for helpful discussions.

\section{The coordinates $\zc$ on $N_w$}
\label{sec_cor}

Let $w \in W$ be a Weyl group element with length $l = l(w)$. Set
\[
R^{+}_{w} ~ = ~ \{ \alpha > 0: ~~ w^{-1} \alpha < 0 \} 
\]
and
\[
\fn_w ~ = ~ \text{span}_{\Bbb C} \{E_{\alpha}: ~~ \alpha \in 
R^{+}_{w} \},
\]
and let $N_w$ be the subgroup of $G$ with Lie algebra $\fn_w$. 
Then $N_w = N \cap w N_{-} w^{-1}$, where $N_{-}$ is the ``opposite"
of $N$, and the map
\begin{equation}
\label{eq_jw}
j_w: ~~ N_w \lrw \Sigma_w: ~~ n \Map nw/B
\end{equation}
is a holomorphic diffeomorphism. Thus, any complex coordinate system on
$N_w$ will  give one on $\Sigma_w$. The coordinates 
$\zc$
that we will introduce in this section,
however, will not be complex for a general $w$ with $l(w) > 1$ 
(see Example \ref{exam_sl3}). 
To obtain these coordinates, we make
use of the 
compact real form $K$ of $G$ and the isomorphism
of $G/B$ and $K/T$. They are motivated by the Bruhat-Poisson structure
on $K/T$.

\bigskip
Let $\dw$ be a representative of $w$ in $K$. 
 For $n \in N_w$, let $a_w(n)$ be the $A$-component
of the element $\dw^{-1} n \dw$ for the Iwasawa decomposition,
i.e.,
\begin{equation}
\label{eq_aw}
\dw^{-1} n \dw ~ = ~ k_1 a_w(n) m_1
\end{equation}
for some $k_1 \in K$ and $m_1 \in N$.
Notice that the map $a_w: N_w \rightarrow A$ depends only on $w$
and not on the choice of $\dw$.  

\bigskip
In this section, we will write down the map 
$a_w: N_w \rightarrow A$ explicitly in the 
coordinates $\zc$ (Theorem \ref{thm_a}).
We will also write down the Haar measure of $N_w$ in these
coordinates (Theorem \ref{thm_haar}). As an application, we show how 
Harish-Chandra's formula for the $c$-function also
follows easily as a product of $1$-dimensional integrals.

\bigskip
We now describe the coordinates $\zc$.

Let again $\dw$ be a representative of $w$ in $K$. Let
\[
C_{\dot{w}} ~ = ~ N_w \circ \dw
\]
be the $N_w$-orbit in $K$ through $\dw$ of the action
$(n, k) \rightarrow n \circ k$. Recall that $n \circ k = k_1$
if $nk = k_1 b$ is the Iwasawa decomposition of $nk$ with
$k_1 \in K$ and $b \in AN$. The map 
\[
J_{\dot{w}}: ~~ N_w \lrw C_{\dot{w}}: ~~ n \Map n \circ \dw
\]
is then a diffeomorphism. Moreover, $J_{\dot{w}}$ followed by the
projection map from $C_{\dot{w}} \subset K \subset G$ to $G/B$ is
just the map $j_w$ in (\ref{eq_jw}). Thus $C_{\dot{w}}$ is a lift
of $\Sigma_w$ in $K$.

\bigskip
{\bf The case of $l(w) = 1$.} 
When $w = \sigma_{\gamma}$ is the simple reflection corresponding to
a simple root $\gamma$, we use, for notational simplicity, $N_{\gamma}$
to denote $N_{\sigma_{\gamma}}$ and $\dot{\gamma}$ for 
$\dot{\sigma}_{\gamma}$. In this case, set
\[
\check{H}_{\gamma} ~ = ~ {\frac{2}{\ll \gamma, \gamma \gg}}
H_{\gamma}, \hspace{.5in}
\check{E}_{\gamma} ~ = ~ \sqrt{{\frac{2}{\ll \gamma, \gamma \gg}}}
~ E_{\gamma}, \hspace{.5in}
\check{E}_{-\gamma} ~ = ~ \sqrt{{\frac{2}{\ll \gamma, \gamma \gg}}}
~ E_{-\gamma}.
\]
Then the map
\begin{equation}
\label{eq_psi-gamma}
\psi_\gamma: ~ sl(2, {\Bbb C}) \lrw \fg: ~~
\left( \begin{array}{cc} 1 & 0\\ 0 & -1 \end{array} \right) \Map 
\check{H}_{\gamma}, \hspace{.1in}
\left( \begin{array}{cc} 0 & 1\\ 0 & 0 \end{array} \right) \Map 
\check{E}_{\gamma}, \hspace{.1in}
\left( \begin{array}{cc} 0 & 0\\ 1 & 0 \end{array} \right) \Map 
\check{E}_{-\gamma}
\end{equation}
is a Lie algebra homomorphism. It induces a Lie group homomorphism
\[
\Psi_{\gamma}: ~ SL(2, {\Bbb C}) \lrw G,
\]
and $\Psi_{\gamma}(SU(2)) \subset K$.
For $z \in {\Bbb C}$, let
\[
n_z ~ = ~ \Psi_{\gamma} \left(
\begin{array}{cc} 1 & z \\ 0 & 1 \end{array} \right) ~ = ~ 
\exp(z \check{E}_{\gamma}) ~ \in ~ N_{\gamma}.
\]
The map
\[
{\Bbb C} \lrw N_{\gamma}: ~~ z \Map n_z
\]
is clearly a parametrization of $N_{\gamma}$ by ${\Bbb C}$. This is a
complex coordinate system on $N_{\gamma}$.

 Notice that the
element $\left( \begin{array}{cc} 0 & i \vspace{-.05in}\\ i & 0 \end{array} 
\right)$ in $SU(2)$ is a representative of the non-trivial 
element in the Weyl group of $SL(2, {\Bbb C})$. Let
\[
\dot{\gamma} ~= ~ \Psi_{\gamma} \left( \begin{array}{cc} 0 & i 
\\ i & 0 \end{array} \right) ~ \in ~ K.
\]
It is a representative of $\sigma_{\gamma}$ in $K$. The following
Iwasawa decomposition in $SL(2, {\Bbb C})$ 
\[
\left(
\begin{array}{cc} 1 & z \\ 0 & 1 \end{array} \right)
\left( \begin{array}{cc} 0 & i 
\\ i & 0 \end{array} \right) ~ = ~ 
\left( \begin{array}{cc} 
{\dfrac{i z}{\sqrt{1 + |z|^2}}} & {\dfrac{i}{\sqrt{1 + |z|^2}}}
\vspace{.1in}\\
{\dfrac{i}{\sqrt{1 + |z|^2}}} & {\dfrac{-i \bar{z}}{\sqrt{1 + |z|^2}}} 
\end{array} \right)  
\left( \begin{array}{cc} \sqrt{1 + |z|^2} & 
{\dfrac{ \bar{z}}{\sqrt{1 + |z|^2}}}\vspace{.1in} \\
\smallskip
0 & {\dfrac{1}{\sqrt{1 + |z|^2}}}
\end{array} \right)
\]
gives the Iwasawa decomposition of $n_z \dot{\gamma}$ in $G$ as
\[
n_z \dot{\gamma} ~ = ~ \Psi_{\gamma}
\left( \begin{array}{cc} 
{\dfrac{i z}{\sqrt{1 + |z|^2}}} & {\dfrac{i}{\sqrt{1 + |z|^2}}} 
\vspace{.1in}\\
{\dfrac{i}{\sqrt{1 + |z|^2}}} & {\dfrac{-i \bar{z}}{\sqrt{1 + |z|^2}}} 
\end{array} \right)  
\Psi_{\gamma} \left( \begin{array}{cc} \sqrt{1 + |z|^2} & 
{\dfrac{ \bar{z}}{\sqrt{1 + |z|^2}}} \vspace{.1in}\\
0 & {\dfrac{1}{\sqrt{1 + |z|^2}}}
\end{array} \right).
\]
Thus
\[
n_z \circ \dot{\gamma} ~ = ~ \Psi_{\gamma}
\left( \begin{array}{cc} 
{\dfrac{i z}{\sqrt{1 + |z|^2}}} & {\dfrac{i}{\sqrt{1 + |z|^2}}}
\vspace{.1in} \\
{\dfrac{i}{\sqrt{1 + |z|^2}}} & {\dfrac{-i \bar{z}}{\sqrt{1 + |z|^2}}} 
\end{array} \right) ~ \in ~ K.
\]
The map 
\[
{\Bbb C} \lrw C_{\dot{\gamma}}:~~ z \Map n_z
\circ \dot{\gamma} ~ = ~ 
\Psi_{\gamma}
\left( \begin{array}{cc}
{\dfrac{i z}{\sqrt{1 + |z|^2}}} & {\dfrac{i}{\sqrt{1 + |z|^2}}}
\vspace{.1in} \\
{\dfrac{i}{\sqrt{1 + |z|^2}}} & {\dfrac{-i \bar{z}}{\sqrt{1 + |z|^2}}}
\end{array} \right)
\]
is a parametrization of $C_{\dot{\gamma}}$ by
$\{z, ~ \bar{z}\}$.
We also see that
\[
a_{\sigma_{\gamma}} (n_z) ~ = ~ \exp({\dfrac{1}{2}} \log(1+|z|^2) 
\check{H}_{\gamma}).
\]

\bigskip
{\bf The general case.} For a general element $w \in W$, let
$l = l(w)$ be the length of $w$, and let
\begin{equation}
\label{eq_rew}
w ~ = ~ \sigma_{\gamma_1} \sigma_{\gamma_2} \cdots 
\sigma_{\gamma_{l}}
\end{equation}
be a reduced decomposition, where $\gamma_1, \gamma_2, ..., 
\gamma_{l}$ are simple roots. Again for notational 
simplicity, we use $N_{\gamma_j}$ to denote 
$N_{\sigma_{\gamma_j}}$ for $j = 1, ..., l$.
We now have the Lie group homomorphism
\[
\Psi_{\gamma_j}: ~~ SL(2, {\Bbb C}) \lrw G
\]
for each $j$. Let again
\[
\dot{\gamma}_j ~ = ~ \Psi_{\gamma_j} 
\left( \begin{array}{cc} 0 & i \\ i & 0 \end{array} \right) ~ \in~ K.
\]
Write an element in $N_{\gamma_j}$ as
\[
n_{z_j} ~ = ~ \exp(z_j \check{E}_{\gamma_j})
\]
for $z_j \in {\Bbb C}$. Set
\[
\dot{w} ~ = ~ \dot{\gamma_1} \dot{\gamma_2} \cdots \dot{\gamma_{\sl}} \in K.
\]

\begin{thm}
\label{thm_cor}
There is a diffeomorphism (between real manifolds)
\[
F_w: ~~ N_{\gamma_1} \times N_{\gamma_2} 
\times \cdots \times N_{\gamma_{\sl}} \lrw N_w
\]
characterized by
\begin{equation}
\label{eq_F}
F_w(n_1, n_2, ..., n_{\sl}) \circ \dot{w} ~ = ~ 
(n_1 \circ \dot{\gamma_1}) 
(n_2 \circ \dot{\gamma_2}) \cdots  
(n_{\sl} \circ \dot{\gamma_{\sl}}) ~ \in ~ K.
\end{equation}
The map
\begin{equation}
\label{eq_F2}
{\Bbb C}^{l} \lrw N_w: ~~ (z_1,~  z_2,~  ...,~ z_{\sl}) \Map
F_w(n_{z_1}, ~ n_{z_2},~  ...~ ,  ~n_{z_{\sl}})
\end{equation}
gives coordinates $\zc$ 
on $N_w$ 
(as a real manifold).
\end{thm}
 
\begin{rem}
\label{rem_soibel}
{\em
This is the same as saying that the map
\[
C_{\dot{\gamma_1}} \times C_{\dot{\gamma_2}} 
\times \cdots \times C_{\dot{\gamma_{\sl}}} \lrw
C_{\dot{w}}
\]
given by multiplication in $K$ is a diffeomorphism. 
This statement is given in \cite{soi:compact}. The 
proof we give below contains a recursive formula for
$F_w$ that will be used later.
}
\end{rem}

\begin{rem}
\label{rem_complex}
{\em
We use $\zc$ instead of $\{z_1, ~ z_2, ~ ..., ~ z_{\sl}\}$
to denote the coordinates to emphasize the fact that 
the map in (\ref{eq_F2}) is in general not a holomorphic diffeomorphism.
See Example \ref{exam_sl3}.
}
\end{rem}

\begin{rem}
\label{rem_choice}
{\em 
Even though we use $F_w$ to denote the map in Theorem \ref{thm_cor},
it depends not only on $w$ but also on the reduced 
decomposition $w = \sigma_{\gamma_1} \sigma_{\gamma_2} \cdots
\sigma_{\gamma_{l}}$ for $w$. Therefore, the coordinates
$\zc$ depend on the choice of the reduced decomposition.
}
\end{rem}

\noindent
{\bf Proof of Theorem \ref{thm_cor}.} We need to show that for  every point
$(n_1, ..., n_{\sl}) \in N_{\gamma_1} \times \cdots
\times N_{\gamma_{\sl}}$, there exists (a necessarily unique)
$F(n_1, ..., n_{\sl}) \in N_w$ such that (\ref{eq_F})
is
satisfied. We also need to show that each $n \in N_w$ 
arises this way.  We prove this by induction on $l(w)$.
When $l(w) = 1$,
the map $F_w$ is the identity map. Now for $w$ with
$l(w) >1$ and with the reduced decomposition given in
(\ref{eq_rew}), set
\[
w_1 ~ = ~ \sigma_{\gamma_1} \sigma_{\gamma_2} \cdots
\sigma_{\gamma_{l-1}}
\]
so that $w = w_1 \sigma_{\gamma_{\sl}}$. We wish to relate
$F_w$ and $F_{w_1}$. To this end, we 
first recall that $N_{w_1} \subset N_w$ \cite{j:group}.
In fact, the multiplication map in $N$ gives a diffeomorphism
\[
N_{w_1} \times \dw_1 N_{\gamma_{\sl}} \dw_{1}^{-1} \lrw N_w,
\]
where  
\[
\dw_1 ~ = ~ \dot{\gamma_1} \dot{\gamma_2} \cdots \dot{\gamma}_{\sl-1}.
\]
Now given $(n_1, n_2, ..., n_{\sl}) \in
N_{\gamma_1} \times N_{\gamma_2} \times \cdots \times
N_{\gamma_{\sl}}$, let
\[
n' ~ = ~ F_{w_1}(n_1, n_2, ..., n_{\sl -1}) \in N_{w_1}
\]
be such that
\[
n' \circ (\dot{\gamma_1} \dot{\gamma_2} \cdots \dot{\gamma}_{l-1}) 
~ = ~ 
(n_1 \circ \dot{\gamma_1})
(n_2 \circ \dot{\gamma_2}) \cdots
(n_{\sl -1 } \circ \dot{\gamma}_{\sl-1}) ~ \in ~ K.
\]
We search for $x \in \dw_1 N_{\gamma_{\sl}} \dw_{1}^{-1}$
such that $n = n' x \in N_w$
satisfies
\begin{equation}
\label{eq_n}
n \circ (\dot{\gamma_1} \dot{\gamma_2} \cdots \dot{\gamma_{\sl}}) ~ = ~ 
(n_1 \circ \dot{\gamma_1})
(n_2 \circ \dot{\gamma_2}) \cdots
(n_{\sl} \circ \dot{\gamma_{\sl}}) ~ \in ~ K.
\end{equation}
Now
\beqa
n \dot{\gamma}_1 \cdots \dg & = & n' x ~\dot{\gamma}_1 \cdots ~
\dgm \dg \\
& = & (n' \dot{\gamma}_1 \cdots \dgm) ~ (\dw_{1}^{-1} x \dw_1)~ \dg \\
& = & (n' \dot{\gamma}_1 \cdots \dgm) ~x' ~\dg,
\eeqa
where
\[
x' ~ = ~ \dw_{1}^{-1} x \dw_1 \in N_{\gamma_{\sl}}.
\]
In the notation we have introduced so far, we have
\begin{equation}
\label{eq_n-prime}
n' \dot{\gamma}_1 \cdots \dgm ~ = ~ (n_1 \circ \dot{\gamma_1})
(n_2 \circ \dot{\gamma_2}) \cdots
(n_{\sl-1} \circ \dot{\gamma}_{\sl-1}) a_{w_1}(n') m'
\end{equation}
for some $m' \in N$.
Thus
\[
n \dot{\gamma}_1 \cdots \dg ~ = ~ \sdgm a_{w_1}(n') ~m' ~x' ~\dg.
\]
Denote by $N_{\hat{\gamma}_{\sl}}$ the subgroup of $N$ with Lie algebra
\[
{\frak n}_{\hat{\gamma}_{\sl}} ~ = ~ \text{span}_{\Bbb C} 
\{E_{\alpha}: ~ \alpha > 0, ~ \alpha \neq \gamma_{\sl} \}.
\]
Then the multiplication map in $N$ induces a diffeomorphism
$N_{\gamma_{\sl}} \times \Nhat \rightarrow N$, so we have the decomposition 
$N = N_{\gamma_{\sl}} \Nhat$. Moreover, $\dg^{-1} \Nhat \dg = 
\Nhat$. Thus, if we take $x'$ to be the element $n^{\gamma_{\sl}} 
\in N_{\gamma_{\sl}}$ in the decomposition
\begin{equation}
\label{eq_hat}
(m')^{-1} a_{w_1}(n')^{-1} n_{\sl} a_{w_1}(n') ~ = ~ 
n^{\gamma_{\sl}} \nhat  \in N_{\gamma_{\sl}} \Nhat,
\end{equation}
then
\[
m' x' ~ = ~ a_{w_1}(n')^{-1} n_{\sl} a_{w_1}(n') (\nhat)^{-1},
\]
and 
\beqa
n \dot{\gamma}_1 \cdots \dg & = & \sdgm ~ n_{\sl} ~  a_{w_1}(n')
(\nhat)^{-1} ~ \dg \\
& = & \sdgm n_{\sl} \dg \left( \dg^{-1} a_{w_1}(n') \dg \right) 
\left(\dg^{-1} (\nhat)^{-1} \dg \right).
\eeqa
Set
\[
m'' ~ = ~ \dg^{-1} (\nhat)^{-1} \dg \in N_{\hat{\gamma}}
\]
and let
\[
n_{\sl} \dg ~ = ~ (n_{\sl} \circ \dg) ~a_{\sl} m_{\sl}
\]
be the Iwasawa decomposition for $n_{\sl} \dg$, so $a_{\sl} \in A$
and $m_{\sl} \in N_{\gamma_{\sl}}$. Then
\[
n \dot{\gamma}_1 \cdots \dg ~ = ~ \sdg ~ a_{\sl} ~ m_{\sl} 
\left(\dg^{-1} a_{w_1}(n') \dg \right) m'',
\]
or
\begin{equation}
\label{eq_n-iwa}
n \dot{\gamma}_1 \cdots \dg ~ = ~ \sdg  ~ a_{\sl} ~
(\dg^{-1} a_{w_1}(n')  \dg ) ~  m''',
\end{equation}
where
\[
m''' ~ = ~ \left( \dg^{-1} a_{w_1}(n') \dg \right)^{-1}
m_{\sl} \left( \dg^{-1} a_{w_1}(n')  \dg \right) m'' \in 
N_{\gamma_{\sl}} \Nhat ~ = ~ N.
\]
Therefore, with this choice of $x' \in N_{\gamma_{\sl}}$, the
element $n = n' x = n' \dw_1 x' \dw_{1}^{-1} \in N_w$
satisfies (\ref{eq_n}). 
Notice that since $N_{\gamma_{\sl}}$ normalizes
$\Nhat$, we can first decompose $(m')^{-1} \in N$ with respect to
the decomposition $N = N_{\gamma_{\sl}} \Nhat$ to get
\begin{equation}
\label{eq_m-n}
(m')^{-1} ~ = ~ m^{\gamma_{\sl}} m^{\hat{\gamma}_{\sl}}
\end{equation}
with $m^{\gamma_{\sl}} \in N_{\gamma}$ and $m^{\hat{\gamma}_{\sl}}
\in \Nhat$. Then 
\begin{equation}
\label{eq_x-prime}
x' ~ = ~ \nhat ~ = ~  m^{\gamma_{\sl}} 
a_{w_1}(n')^{-1} n_{\sl} a_{w_1}(n') \in N_{\gamma_{\sl}}.
\end{equation}
To summerize, we set
\begin{equation}
\label{eq_f}
f: ~~ N_{w_1} \times N_{\gamma_{\sl}} \lrw N_w: ~~ (n', ~ n_{\sl})
\Map n' \dw_1 (m^{\gamma_{\sl}}
a_{w_1}(n')^{-1} n_{\sl} a_{w_1}(n')) \dw^{-1}_{1} ~ \in ~ N_{w},
\end{equation}
and define
\begin{equation}
\label{eq_n-final}
F_w (n_1, ..., n_{\sl}) ~  =  ~  n
~ = ~ n' \dw_1 (m^{\gamma_{\sl}}
a_{w_1}(n')^{-1} n_{\sl} a_{w_1}(n')) \dw^{-1}_{1} \in N_w.
\end{equation}
Then $F_w: N_{\gamma_1} \times N_{\gamma_2} \times \cdots
\times N_{\gamma_{\sl}} \rightarrow N_w$ is a well-defined map and 
\[
F_w ~ = ~ f \circ (F_{w_1} \times id).
\]
It is also clear now that $F_w$ is a diffeomorphism.
\qed

Formula (\ref{eq_n-iwa}) also gives the following recursive formula
for $a_w(n)$:
\begin{equation}
\label{eq_a-rec}
a_w (n) ~ = ~ a_{\sl} ~ \dg^{-1} a_{w_1}(n') \dg,
\end{equation}
where $n_{\sl} \dg = (n_{\sl} \circ \dg) a_{\sl} m_{\sl}$ is
the Iwasawa decomposition for $n_{\sl} \dg$. Now for 
each $j = 1, ..., l$, let 
\[
n_j \dot{\gamma}_j ~ = ~ (n_j \circ \dot{\gamma}_j) ~ a_j m_j
\]
be the Iwasawa decomposition. We get from (\ref{eq_a-rec}) that
\begin{equation}
\label{eq_a-exp}
a_w(n) ~ = ~ a_{\sl} (\dg^{-1} a_{\sl -1} \dg) 
(\dg^{-1} \dot{\gamma}_{{\sl}-1}^{-1} a_{\sl -2} \dgm \dg) \cdots
(\dg^{-1} \dot{\gamma}_{{\sl}-1}^{-1} \cdots \dot{\gamma}_{2}^{-1} a_1 
\dot{\gamma}_2 \cdots \dgm \dg),
\end{equation}
or, since $A$ is commutative,
\[
\dw a_w(n) \dw^{-1}  ~ = ~  \prod_{j=1}^{l}
\sigma_{\gamma_1} \sigma_{\gamma_2} \cdots \sigma_{\gamma_{j-1}}
\sigma_{\gamma_j} (a_j).
\]
We know that in the $\zc$ coordinates, $a_j$ is given by
\[
a_j ~ = ~ \exp({\frac{1}{2}} (1 + |z_j|^2) \check{H}_{\gamma_j})
\]
for $j = 1, ..., l$. Thus
\beqa
\sigma_{\gamma_1} \sigma_{\gamma_2} \cdots \sigma_{\gamma_{j-1}}
\sigma_{\gamma_j} (a_j) & = & 
\sigma_{\gamma_2} \cdots \sigma_{\gamma_{j-1}} \exp
(-{\frac{1}{2}} \log(1+|z_j|^2) \check{H}_{\gamma_j})\\
& = & \exp (-{\frac{1}{2}} \log(1+|z_j|^2) \check{H}_{\alpha_j}),
\eeqa
where, for each $j = 1, ..., l$,
\begin{equation}
\label{eq_alphaj}
\alpha_j ~ = ~ \sigma_{\gamma_1} 
\sigma_{\gamma_2} \cdots \sigma_{\gamma_{j-1}} (\gamma_j).
\end{equation}
Recall that 
\[
\{ \alpha_1, ~ \alpha_2, ~ \cdots ~ \alpha_{\sl} \} ~ = ~ 
R^{+}_{w} ~ = ~ \{ \alpha > 0: ~~ w^{-1} \alpha < 0 \}.
\]
Thus we have
\begin{equation}
\label{eq_a-w}
\dw a_w(n) \dw^{-1} ~ = ~ \prod_{j=1}^{l} 
\exp (-{\frac{1}{2}} \log(1+|z_j|^2) \check{H}_{\alpha_j}).
\end{equation}
Let 
\begin{equation}
\label{eq_beta}
\beta_j ~ = ~ - w^{-1} \alpha_j ~ = ~ -\sigma_{\sl} \sigma_{\sl -1}
\cdots \sigma_j(\gamma_j) ~ = ~ \sigma_{\sl} \sigma_{\sl -1}
\cdots \sigma_{j+1}(\gamma_j),
\end{equation}
i.e.,
\beqa
& & \beta_1 ~ = ~ \sigma_{\gamma_{\sl}} \sigma_{\gamma_{\sl -1}} 
\cdots \sigma_{\gamma_2} (\gamma_1) \\
& & \beta_2 ~ = ~ \sigma_{\gamma_{\sl}} \sigma_{\gamma_{\sl -1}} 
\cdots \sigma_{\gamma_3} (\gamma_2) \\ 
& & \hspace{.3in} \cdots \\
& & \beta_{\sl -1} ~ = ~ \sigma_{\gamma_{\sl}} (\gamma_{\sl -1}) \\
& & \beta_{\sl} ~ = ~ \dg.
\eeqa
We then know that
\[
\{\beta_1, ~ \beta_2, ~ ...~ , \beta_{\sl} \} ~ = ~ R_{w^{-1}}^{+} 
~ = ~ \{\beta > 0: ~~ w \beta < 0\}.
\]
We also know that $\ll \beta_j, ~ \beta_j \gg = \ll \alpha_j, ~ \alpha_j \gg
= \ll \gamma_j, ~ \gamma_j \gg$ for each $j = 1, ..., l$. This fact will
be used later.

\bigskip
The following theorem now follows immediately from 
(\ref{eq_a-w}).

\begin{thm}
\label{thm_a}
In the $\zc$-coordinates, the function $a_w$ on $N_w$ defined by
(\ref{eq_aw}) is explicitly given by
\begin{equation}
\label{eq_aw-zc}
a_w(n) ~ = ~ \prod_{j=1}^{l} \exp({\frac{1}{2}} \log(1 + |z_j|^2) 
\check{H}_{\beta_j}),
\end{equation}
where the $\beta_j$'s are given by (\ref{eq_beta}).
\end{thm}

\begin{rem}
\label{rem_pickrell}
{\em
This formula for $a_w(n)$ also follows from a product formula
found by Doug Pickrell in \cite{pi:schubert}.
}
\end{rem}

We now look at the left invariant Haar measure on $N_w$.

\begin{thm}
\label{thm_haar}
In the $\zc$-coordinates, a left invariant (and thus also bi-invariant)
Haar measure on $N_w$ is given as
\begin{equation}
\label{eq_dn}
dn ~ = ~ \lambda_w \prod_{j=1}^{l} 
 ~(1 + |z_j|^2)^{{\dfrac{2 \ll \rho, ~ \beta_j \gg}
{\ll \beta_j, ~ \beta_j \gg}} -1} dz_j \wedge d \bar{z}_j,
\end{equation}
where 
\begin{equation}
\label{eq_lambda}
\lambda_w ~ = ~ {\frac{1}{(-2\pi i)^l}} \prod_{j=1}^{l} 
{\frac{2 \ll \rho, ~ \beta_j \gg}{\ll \beta_j, ~ \beta_j \gg }}
~ = ~ \prod_{j=1}^{l} {\frac{i \ll \rho, ~ \beta_j \gg}{ \pi
\ll \beta_j, ~ \beta_j \gg}}
\end{equation}
is such that
\[
\int_{N_w} a_w(n)^{-4 \rho} dn ~ = ~ 1.
\]
\end{thm}
 
\noindent
{\bf Proof.} Again we prove by induction on $l(w)$. When
$l(w) = 1$ so $w = \sigma_{\gamma}$ for a simple root $\gamma$, 
we have
\[
d n_{\gamma} ~ = ~ -{\frac{1}{2 \pi i}} d z_1 \wedge d \bar{z}_1.
\]
Now for $w$ with $l(w) = l >1$, we use the same notation as that in
the proof of Theorem \ref{thm_cor}. Let
\[
N_{\alpha_{\sl}} ~ = ~ \dw_1 N_{\gamma_{\sl}} \dw_{1}^{-1}
\]
be the subgroup of $N$ with Lie algebra ${\Bbb C} E_{\alpha_{\sl}}$,
where $\alpha_{\sl} = w_1 (\gamma_{\sl})$. Then the 
multiplication map
\[
\mu: ~~ N_{w_1} \times N_{\alpha_{\sl}} \lrw N: ~~ (n', ~ 
n_{\alpha_{\sl}}) \Map n' n_{\alpha_{\sl}}
\]
is a diffeomorphism. Since $N$ is unipotent,
we have, under the map $\mu$,
\[
dn ~ = ~ \lambda ~  dn' d n_{\alpha_{\sl}},
\]
where $\lambda$ is a constant to be determined later, and we take
\[
d n_{\alpha_{\sl}} ~ = ~ du \wedge d \bar{u}
\]
if $N_{\alpha_{\sl}}$ is parametrized by $\{u, ~ \bar{u}\}$ via 
$n_{\alpha_{\sl}} = \exp(u E_{\alpha_{\sl}})$.

Consider now the parametrization of $N_w$ by ${\Bbb C}^l$ via
$F_w$. Write $n = F_w(z_1, \bar{z}_1, ..., 
z_{\sl}, \bar{z}_{\sl})$ if $n = F_w(n_1, ..., n_{\sl})$
where $n_j = \exp(z_j \check{E}_{\gamma_j})$ for each $j$.
Let again 
\[
n' ~ = ~ F_{w_1} (z_1, ~ \bar{z}_1, ~ ...~ , z_{\sl -1}, ~ 
\bar{z}_{\sl -1}) \in N_{w_1}.
\]
Recall that the element $m' \in N$ is given in (\ref{eq_n-prime}) and that
$m^{\gamma_{\sl}} \in N_{\gamma_{\sl}}$ is given in (\ref{eq_m-n}). 
Write 
\[
m^{\gamma_{\sl}} ~ = ~ \exp(m(z_1, \bar{z}_1, ..., z_{\sl -1},
\bar{z}_{\sl -1}) \check{E}_{\gamma_{\sl}}).
\]
Then we know from (\ref{eq_n-final}) that 
\[
n ~ = ~ n' \exp\left( ( m(z_1, \bar{z}_1, ..., z_{\sl -1},
\bar{z}_{\sl -1}) + a_{w_1}(n')^{-\gamma_{\sl}} z_{\sl}) 
 v_{\sl} {E}_{\alpha_{\sl}} \right),
\]
where $v_{\sl} \in {\Bbb C}$ is such that $Ad_{\dot{w}_1} 
\check{E}_{\gamma_{\sl}} = v_{\sl} E_{\alpha_{\sl}}$.
Set 
\[
u ~ = ~v_{\sl} ( m(z_1, \bar{z}_1, ..., z_{\sl -1},
\bar{z}_{\sl -1}) + a_{w_1}(n')^{-\gamma_{\sl}} z_{\sl}).
\]
Assume that $d n'$ is given as in the theorem for $w_1$. This means,
noting the definition of the $\beta_j$'s, that
\beqa
d n' & = & \lambda_{w_1} \prod_{j=1}^{l -1}~
(1 + |z_j|^2)^{{\dfrac{2 \ll \rho, ~ \sigma_{\gamma_{\sl}} \beta_j \gg}
{\ll \sigma_{\gamma_{\sl}} \beta_j, ~ \sigma_{\gamma_{\sl}} \beta_j \gg}} 
~ -~  1} dz_j \wedge d \bar{z}_j \\
& = & \lambda_{w_1} \prod_{j=1}^{l -1} 
~(1 + |z_j|^2)^{- {\dfrac{2 \ll \gamma_{\sl}, ~ \beta_j \gg}{ 
\ll \beta_j, ~ \beta_j \gg}}} 
(1 + |z_j|^2)^{{\dfrac{2 \ll \rho, ~ \beta_j \gg}{ \ll \beta_j, ~ 
\beta_j \gg}}~ - ~ 1} dz_j \wedge d \bar{z}_j.
\eeqa
Here we have just used the fact that $\rho - \sigma_{\gamma_{\sl}}
\rho = \gamma_{\sl}$. On the other hand, by Theorem
\ref{thm_a}, we have
\[
a_{w_1}(n') ~ = ~ \prod_{j=1}^{l -1} \exp({\frac{1}{2}} \log(1 + |z_j|^2)
\check{H}_{\sigma_{\gamma_{\sl}}(\beta_j)}),
\]
so
\[
a_{w_1}(n')^{-2 \gamma_{\sl}} ~ = ~ \prod_{j=1}^{l -1}
~(1+|z_j|^2)^{ 
{\dfrac{\ll -2 \gamma_{\sl}, ~  \sigma_{\gamma_{\sl}} \beta_j \gg}
{\ll \sigma_{\gamma_{\sl}} \beta_j, ~ 
\sigma_{\gamma_{\sl}} \beta_j \gg}} } ~ = ~
 \prod_{j=1}^{l -1}~  (1 + |z_j|^2)^{
{\dfrac{2 \ll \gamma_{\sl}, ~ \beta_j \gg }{\ll \beta_j, ~ \beta_j \gg}}}.
\]
Therefore,
\beqa
dn & = & \lambda ~ dn' ~ d n_{\alpha_{\sl}}
~ = ~ \lambda'~ a_{w_1}(n')^{-2 \gamma_{\sl}} ~  d n' ~(dz_{\sl} \wedge
d \bar{z}_{\sl}) \\
& = & \lambda' \lambda_{w_1} \left( \prod_{j=1}^{l -1}
~(1 + |z_j|^2)^{{\dfrac{2 \ll \rho, ~ \beta_j \gg}{ \ll \beta_j, ~
\beta_j \gg}}~ - ~ 1} dz_j \wedge d \bar{z}_j \right) \wedge 
(dz_{\sl} \wedge d \bar{z}_{\sl}),
\eeqa
where $\lambda' = \lambda |v_{\sl}|^2$ is a new constant to be determined 
later. Since $\beta_{\sl} = \gamma_{\sl}$ is a simple root, 
we have
\[
{\frac{2\ll \rho, ~ \gamma_{\sl} \gg}{\ll \gamma_{\sl}, ~ \gamma_{\sl} \gg}}
~ = ~ 1.
\]
Thus 
\[
dn ~ = ~ \lambda' \lambda_{w_1} 
\prod_{j=1}^{l}
 ~(1 + |z_j|^2)^{{\dfrac{2 \ll \rho, ~ \beta_j \gg}
{\ll \beta_j, ~ \beta_j \gg}} ~ - ~ 1} dz_j \wedge d \bar{z}_j.
\]
Using Theorem \ref{thm_a}, we see that the integral 
\[
\int_N a_w(n)^{-4 \rho} dn
\]
is now a product of $1$-dimensional ones and is easily calculated. 
The constant $\lambda_w = \lambda' \lambda_{w_1}$ must be 
given by (\ref{eq_lambda}) for the above integral to be
equal to $1$.
\qed

\begin{exam}
\label{exam_c}
{\em We recall that for the complex group $G$ considered as
a real Lie group, the {\bf $c$-function} for a Weyl group element
$w$ is defined to be (see \cite{h:gga}, Chapter IV, $\S 6$)
\[
c_w(\lambda) ~ = ~ \int_{\bar{N}_w} a(\bar{n})^{-(i \lambda + 2 \rho)}
d \bar{n},
\]
where $\bar{N}_w \subset N_{-}$ is the ``opposite" of $N_w$ and
$a(\bar{n})$ is the $A$-component in the Iwasawa decomposition of $\bar{n}
\in \bar{N}_w$. In our notation, we have
\[
c_{w^{-1}} (\lambda) ~ = ~ \int_{N_w} a_w(n)^{-(i \lambda + 2 \rho)}
dn.
\]
Using our formulas for $a_w(n)$ and for $dn$ in the $\zc$ coordinates,
one immediately reduces the integral to a product of
$1$-dimensional ones and  gets
\beqa
c_{w^{-1}} (\lambda) & = & \lambda_{w} \prod_{j=1}^{l(w)}
\int_{{\Bbb R}^2} (1 + |z_j|^2)^{- ~ {\dfrac{\ll i \lambda + 2 \rho, ~
\beta_j \gg}
{\ll \beta_j, ~ \beta_j \gg}} ~ + ~ 
{\dfrac{2 \ll \rho, ~ \beta_j \gg}{\ll \beta_j, ~ \beta_j \gg}} ~ - ~ 1}
dz_j \wedge d \bar{z}_j \\
& = & \lambda_{w} \prod_{j=1}^{l(w)} \int_{{\Bbb R}^2}
(1 + x_{j}^{2} + y_{j}^{2})^{- ~ {\dfrac{\ll i \lambda, ~ \beta_j \gg}
{\ll \beta_j, ~ \beta_j \gg}} ~ - ~ 1} (-2i) dx_j dy_j \\
& = & \lambda_{w} \prod_{j=1}^{l(w)}
(~-2 \pi i) \int_{0}^{\infty}
(1+r_{j}^{2})^{- ~ {\dfrac{\ll i \lambda, ~ \beta_j \gg}
{\ll \beta_j, ~ \beta_j \gg}} ~ - ~ 1} 2 r_j d r_j \\
& = & (-2 \pi i)^{l(w)}~ \lambda_w \prod_{j=1}^{l(w)}
{\frac{\ll \beta_j, ~ \beta_j \gg}{ \ll i \lambda, ~ \beta_j \gg}} 
\hspace{.4in}
\text{if} ~ \text{Re} \ll i\lambda, ~ \beta_j \gg ~ > 0 \hspace{.1in}
\text{for each} ~ j\\
& = & \prod_{j=1}^{l(w)} {\frac{\ll 2 \rho, ~ \beta_j \gg }{\ll
i \lambda, ~ \beta_j \gg}} \hspace{.4in}
\text{if} ~ \text{Re}\ll i\lambda, ~ \beta_j \gg ~ > 0 \hspace{.1in}
\text{for each} ~ j.
\eeqa
This is the well-known formula of Harish-Chandra (see Theorem 5.7 in
\cite{h:gga}). As we have mentioned in the Introduction, our calculation here
is easier because we have pushed the induction argument that is normally
used in the calculations for the $c$-functions into the calculations for
$a_w(n)$ and $dn$.
}
\end{exam}

\begin{exam}
\label{exam_sl3}
{\em
Consider the example of $\fg = sl(3, {\Bbb C})$. We take $w$ be to
the longest Weyl group element $w_0 = (1,2)(2,3)(1,2)$. 
In this case, parametrize $N_{w_0} = N$ by complex
coordinates
\[
(u_1, ~ u_2, ~ u_3) \Map n ~ = ~  \left(
\begin{array}{ccc} 1 & u_1 & u_3 \\
                   0 & 1 & u_2 \\
                   0 & 0 & 1 \end{array} \right).
\]
We have
\[
\dot{\gamma}_1 ~ = ~ \left( \begin{array}{ccc} 0 & i & 0 \\ i & 0 & 0\\ 
0 & 0 & 1 \end{array} \right)
\hspace{.2in}
\dot{\gamma}_2 ~ = ~ \left( \begin{array}{ccc} 1 & 0 & 0 \\ 0 & 0 & i\\ 
0 & i & 0 \end{array} \right)
\hspace{.2in}
\dot{\gamma}_3 ~ = ~ \left( \begin{array}{ccc} 0 & i & 0 \\ i & 0 & 0\\ 
0 & 0 & 1 \end{array} \right),
\]
so
\[
\dw  ~ = ~ \dot{\gamma}_1 \dot{\gamma}_2 \dot{\gamma}_3 ~ = ~ 
\left( \begin{array}{ccc} 0 & 0 & -1 \\ 0 & -1 & 0\\ 
-1 & 0 & 0 \end{array} \right).
\]
The Iwasawa decomposition of $n \dw$ in $SL(3, {\Bbb C})$ is 
\beqa
& & n \dw  ~ = ~ 
\left( \begin{array}{ccc} -u_3 & -u_1 & -1 \\
-u_2 & -1 & 0 \\ -1 & 0 & 0 \end{array} \right) ~ = ~ \vspace{.1in} \\
& &  \left( \begin{array}{ccc} 
        -{\dfrac{u_3}{\Delta_3}} & -{\dfrac{u_1(1+|u_2|^2) -\bar{u}_2 u_3}
        {\Delta_2 \Delta_3}} & -{\dfrac{1}{\Delta_2}} \vspace{.1in} \\
        -{\dfrac{u_2}{\Delta_3}} & -{\dfrac{1+|u_2|^2 - u_1 u_2 \bar{u}_3}
        {\Delta_2 \Delta_3}} & {\dfrac{\bar{u}_1}{\Delta_2}} \vspace{.1in}\\
        -{\dfrac{1}{\Delta_3}} & {\dfrac{\bar{u}_2 + u_1 \bar{u}_3}
       {\Delta_2 \Delta_3}} & {\dfrac{\bar{u}_3 - \bar{u}_1 \bar{u}_2}
       {\Delta_2}} \end{array} \right)
\left( \begin{array}{ccc}
\Delta_3 & {\dfrac{\bar{u}_2 + u_1 \bar{u}_3}{\Delta_3}} & 
{\dfrac{\bar{u}_3}{\Delta_3}} \vspace{.1in} \\
0 & {\dfrac{\Delta_2}{\Delta_3}} & {\dfrac{\bar{u}_1(1+|u_2|^2)-u_2 \bar{u}_3}
{\Delta_2 \Delta_3}} \vspace{.1in} \\ 0 & 0 & {\dfrac{1}{\Delta_2}} \end{array}
\right)
\eeqa
where
\[
\Delta_1 ~ = ~ \sqrt{1+|u_1|^2}, \hspace{.3in}
\Delta_2 ~ = ~ \sqrt{1+|u_1|^2 + |u_1 u_2 - u_3|^2}, \hspace{.3in}
\Delta_3 ~ = ~ \sqrt{1+|u_2|^2 + |u_3|^2}. 
\]
On the other hand, for $z_1, z_2$ and $z_3$ in $\Bbb C$, let 
\[
n_1 ~ = ~ \left( \begin{array}{ccc} 1 & z_1 & 0 \\ 0 & 1 & 0 \\
0 & 0 & 1 \end{array} \right) \hspace{.3in}
n_2 ~ = ~ \left( \begin{array}{ccc} 1 & 0 & 0 \\ 0 & 1 & z_2 \\
0 & 0 & 1 \end{array} \right) \hspace{.3in}
n_3 ~ = ~ \left( \begin{array}{ccc} 1 & z_3 & 0 \\ 0 & 1 & 0 \\
0 & 0 & 1 \end{array} \right). 
\]
We have
\beqa
& & (n_1 \circ \dot{\gamma}_1) (n_2 \circ \dot{\gamma}_2) (n_3 \circ 
\dot{\gamma}_3) \vspace{.3in}\\
&  & ~~~~~~ = ~ {\dfrac{1}{\epsilon_1 \epsilon_2 \epsilon_3}}
\left( \begin{array}{ccc} i z_1 & i & 0 \\ i & -i \bar{z}_1 & 0 \\
0 & 0 & \epsilon_1 \end{array} \right) 
\left( \begin{array}{ccc} 
\epsilon_2 & 0 & 0 \\ 0 & i z_2 & i \\ 0 & i & -i \bar{z}_2 
\end{array} \right)
\left( \begin{array}{ccc} i z_3 & i & 0 \\ i & -i \bar{z}_3 & 0 \\
0 & 0 & \epsilon_3 \end{array} \right) \vspace{.3in} \\
&  & ~~~~~~ = ~ {\dfrac{1}{\epsilon_1 \epsilon_2 \epsilon_3}} 
\left( \begin{array}{ccc} 
-\epsilon_2 z_1 z_3 - i z_2 & -\epsilon_2 z_1 + i z_2 \bar{z}_3 & 
- \epsilon_3  \vspace{.1in} \\ -\epsilon_2 z_3 + i \bar{z}_1 z_2 & 
-\epsilon_2 - i \bar{z}_1 z_2 \bar{z}_3 & \bar{z}_1 \epsilon_3 \vspace{.1in}\\
-\epsilon_1 & \epsilon_1 \bar{z}_3 & -i \epsilon_1 \epsilon_3 \bar{z}_2 
\end{array} \right),
\eeqa
where $\epsilon_j = \sqrt{1 + |z_j|^2}$ for $j = 1, 2, 3$.
By setting
\[
n \circ \dw ~ = ~ 
(n_1 \circ \dot{\gamma}_1) (n_2 \circ \dot{\gamma}_2) (n_3 \circ 
\dot{\gamma}_3),
\]
we get the following
coordinate change between our coordinates 
$\{(z_1, \bar{z}_1, z_2, \bar{z}_2, z_3, \bar{z}_3 \}$  and
the $u$'s: 
\smallskip
\[
u_1 ~ = ~ z_1, \hspace{.3in}
u_2 ~ = ~ {\dfrac{\epsilon_2 z_3 - i \bar{z_1} z_2}{\epsilon_1}},
\hspace{.3in}
u_3 ~ = ~ {\dfrac{\epsilon_2 z_1 z_3 + i z_2}{\epsilon_1}},
\]
or
\[
z_1 ~ = ~ u_1, \hspace{.3in}
z_2 ~ = ~ i{\dfrac{u_1 u_2 - u_3}{\Delta_1}}, \hspace{.3in}
z_3 ~ = ~ {\dfrac{\bar{u}_1 u_3 + u_2}{\Delta_2}}.
\]
\vspace{.1in}
Notice that this is not a holomorphic change of coordinates. Thus the 
$\{z_1, \bar{z}_1, z_2, \bar{z}_2, z_3, \bar{z}_3 \}$
coordinates are not complex.
Now in the $u$-coordinates, we have
\[
a_w(n) ~ = ~ \left( \begin{array}{ccc} \Delta_3 & 0 & 0 \vspace{.1in}\\
0 & {\dfrac{\Delta_2}{\Delta_3}} & 0 \vspace{.1in}\\ 0 & 0 & 
{\dfrac{1}{\Delta_2}} \end{array} \right).
\]
Under the coordinate change, we have 
\[
\Delta_1 ~ = ~ \epsilon_1, \hspace{.3in}
\Delta_2 ~ = ~ \epsilon_1 \epsilon_2, \hspace{.3in}
\Delta_3 ~ = ~ \epsilon_2 \epsilon_3.
\]
Thus  
we get $a_w(n)$ in the 
$\{(z_1, \bar{z}_1, z_2, \bar{z}_2, z_3, \bar{z}_3 \}$ coordinates as:
\[
a_w(n) ~= ~ \left( \begin{array}{ccc} \epsilon_2 \epsilon_3 & 0 & 0 
\vspace{.1in}\\
0 & {\dfrac{\epsilon_1}{\epsilon_3}} & 0 \vspace{.1in}\\ 0 & 0 & 
{\dfrac{1}{\epsilon_1 \epsilon_2}} \end{array} \right).
\]
This  is the same as what one would obtain from Theorem \ref{thm_a}. Similarly, 
the left invariant Haar measure $dn$ is, up to a constant multiple,
given in the $u$-coordinates by
\[
dn ~ = ~ du_1 \wedge d \bar{u}_1 \wedge du_2 \wedge
d \bar{u}_2 \wedge du_3 \wedge d \bar{u}_3.
\]
After the change of coordinates, we get
\[
dn ~ = ~ \epsilon_{2}^{2} dz_1 \wedge d \bar{z}_1 \wedge dz_2 \wedge
d \bar{z}_2 \wedge dz_3 \wedge d \bar{z}_3.
\]
Again, this is the same as what one would obtain from Theorem \ref{thm_haar}.
}
\end{exam}

\bigskip
The following proposition will be used in Section \ref{sec_kost-thm}.
\begin{prop}
\label{prop_scalar}
We have
\[
F_w(e, ..., e, n_j, e, ..., e) ~ = ~
(\dot{\gamma}_1 \dot{\gamma}_2 \cdots \dot{\gamma}_{j-1}) ~ n_j ~
(\dot{\gamma}_1 \dot{\gamma}_2 \cdots \dot{\gamma}_{j-1})^{-1}
\in N_{\alpha_j}
\]
for $j = 1, ..., l$ and $n_j \in N_j$, and
\[
(F_{w})_{*}(0) \left({\frac{\partial}{\partial z_1}} \wedge
{\frac{\partial}{\partial \bar{z}_1}} \wedge \cdots \wedge
{\frac{\partial}{\partial z_{\sl}}} \wedge
{\frac{\partial}{\partial \bar{z}_{\sl}}}\right) ~ = ~
(\prod_{j=1}^{l}
{\frac{i}{\ll \alpha_j, ~ \alpha_j \gg}})
E_{\alpha_1} \wedge i E_{\alpha_1} \wedge \cdots \wedge
E_{\alpha_{\sl}} \wedge i E_{\alpha_{\sl}},
\]
where $(F_{w})_{*}(0)$ is the differential of $F_w$ at $0$.
\end{prop}
 
\noindent
{\bf Proof.} Set $n =
(\dot{\gamma}_1 \dot{\gamma}_2 \cdots \dot{\gamma}_{j-1}) n_j
(\dot{\gamma}_1 \dot{\gamma}_2 \cdots \dot{\gamma}_{j-1})^{-1}.$
Let
\[
n_j \dot{\gamma}_j ~ = ~ (n_j \circ \dot{\gamma}_j)
~ a_j ~ m_j
\]
be the Iwasawa decomposition of $n_j \dot{\gamma}_j$, where
$m_j \in N_{\gamma_j}$.
Then
\[
n_j \dot{\gamma}_j \dot{\gamma}_{j+1} \cdots
\dg ~ = ~ (n_j \circ \dot{\gamma}_j)
(\dot{\gamma}_{j+1} \cdots
\dg) (\dot{\gamma}_{j+1} \cdots
\dg)^{-1}a_j ~ (\dot{\gamma}_{j+1} \cdots
\dg ) (\dot{\gamma}_{j+1} \cdots
\dg)^{-1} m_j (\dot{\gamma}_{j+1} \cdots
\dg)
\]
Set
\beqa
a_j' & = & (\dot{\gamma}_{j+1} \cdots
\dg)^{-1} a_j
(\dot{\gamma}_{j+1} \cdots \dg) \in A \\
m_j' & = & (\dot{\gamma}_{j+1} \cdots \dg)^{-1} m_j
(\dot{\gamma}_{j+1} \cdots \dg) \in N.
\eeqa
Then
\[
n ~\dot{\gamma}_1 \cdots \dg  ~ = ~ \dot{\gamma}_1 \cdots
\dot{\gamma}_{j-1}
(n_j \circ \dot{\gamma}_j)~
(\dot{\gamma}_{j+1} \cdots \dg)~  a_j' ~ m_j'.
\]
Thus
\[
n \circ (\dot{\gamma}_1 \cdots \dg) ~ = ~
\dot{\gamma}_1 \cdots
\dot{\gamma}_{j-1}
(n_j \circ \dot{\gamma}_j)
~ \dot{\gamma}_{j+1} \cdots \dg
\]
Hence
\[
F_w(e, ..., e, n_j, e, ..., e) ~ = ~ n.
\]
Write $n_j ~ = ~ \exp(z_j \check{E}_{\gamma_j})$. Then
\beqa
(\dot{\gamma}_1 \dot{\gamma}_2 \cdots \dot{\gamma}_{j-1}) ~ n_j ~
(\dot{\gamma}_1 \dot{\gamma}_2 \cdots \dot{\gamma}_{j-1})^{-1}
& = & \exp(z_j Ad_{\dot{\gamma}_1 \dot{\gamma}_2 \cdots \dot{\gamma}_{j-1}}
(\check{E}_{\gamma_j})) \vspace{.2in} \\
& = & \exp(\sqrt{{\dfrac{2}{\ll \gamma_j, ~ \gamma_j \gg}}} z_j
Ad_{\dot{\gamma}_1 \dot{\gamma}_2 \cdots \dot{\gamma}_{j-1}}
(E_{\gamma_j})).
\eeqa
Since $\sigma_{\gamma_1} \sigma_{\gamma_2} \cdots \sigma_{\gamma_{j-1}}
(\gamma_j) = \alpha_j$ by definition, we have
\[
Ad_{\dot{\gamma}_1 \dot{\gamma}_2 \cdots \dot{\gamma}_{j-1}}
({E}_{\gamma_j}) ~ = ~ c_j E_{\alpha_j}
\]
for some complex number $c_j$. Write $z_j = x_j + i y_j$ and
$c_j = u_j + i v_j$. Using the fact that 
$\ll \gamma_j, ~ \gamma_j \gg =
\ll \alpha_j, ~ \alpha_j \gg$, we get
\beqa
(F_w)_* (0) ({\frac{\partial}{\partial z_j}}
\wedge {\frac{\partial}{\partial \bar{z}_j}}) & = & {\frac{i}{2}}
(F_w)_* (0) ({\frac{\partial}{\partial x_j}} \wedge
{\frac{\partial}{\partial y_j}}) \\
& = & {\frac{i}{\ll \alpha_j, ~ \alpha_j \gg }}
(u_j E_{\alpha_j} + v_j (i E_{\alpha_j}))
\wedge (-v_j E_{\alpha_j} + u_j (i E_{\alpha_j})) \\
& = & {\frac{i}{\ll \alpha_j, ~ \alpha_j \gg }} |c_j|^2
~ E_{\alpha_j} \wedge i E_{\alpha_j}.
\eeqa
Since
both root vectors $E_{\gamma_j}$ and $E_{\alpha_j}$ have length $1$
with respect to the $K$-invariant Hermitian form on $\fg$ induced by
$\fk$,
we have $|c_j|^2 = 1$.  The statement about $(F_w)_*(0)$ now follows
immediately from this.
\qed

We now look at the $T$-action on $N_w$ by conjugations:
\[
T \times N_{w} \lrw N_w: ~~ (t, ~ n) \Map t ~ n ~ t^{-1}.
\]
For a given $t \in T$, set
\beqa
t_1 & = & t \\
t_2 & = & \dot{\gamma}_1 ~ t ~ \dot{\gamma}_{1}^{-1}\\
t_3 & = & \dot{\gamma}_2 \dot{\gamma}_1 ~  t~
(\dot{\gamma}_2 \dot{\gamma}_1)^{-1}\\
& & \hspace{.2in} \cdots \\
t_{\sl} & = & (\dgm \dot{\gamma}_{\sl -2} \cdots \dot{\gamma}_1 )
~ t~  (\dgm \dot{\gamma}_{\sl -2} \cdots \dot{\gamma}_1 )^{-1}.
\eeqa
Equip $N_{\gamma_1} \times N_{\gamma_2} \times \cdots \times
N_{\gamma_{\sl}}$ with the $T$-action given by
\begin{equation}
\label{eq_t-on-n-times}
t \cdot (n_1, ~ n_2, ~ ...,~ n_{\sl}) ~ = ~
(t_1 n_1 t_{1}^{-1}, ~ t_2 n_2 t_{2}^{-1}, ~ \cdots, ~
t_{\sl} n_{\sl} t_{\sl}^{-1}).
\end{equation}
 
\bigskip
\begin{prop}
\label{prop_t-on-N}
1) With respect to the $T$-actions on $N_w$ by conjugation and on
$N_{\gamma_1} \times N_{\gamma_2} \times \cdots \times
N_{\gamma_{\sl}}$ as given by (\ref{eq_t-on-n-times}), the map
$F_w$ is $T$-equivariant;
 
2) In the $\zc$ coordinates, the $T$-action on $N_w$ is given by
\begin{equation}
\label{eq_t-cor}
t \cdot (z_1, ~ \bar{z}_1, ~ ...~ ,  ~ z_{\sl}, ~
\bar{z}_{\sl}) ~ = ~ (t^{\alpha_1} z_1, ~ t^{-\alpha_1} \bar{z}_1, ~
...~ , ~ t^{\alpha_{\sl}} z_{\sl}, ~ t^{-\alpha_{\sl}} \bar{z}_{\sl}),
\end{equation}
where, recall, $\alpha_j = \sigma_{\gamma_1} \sigma_{\gamma_2} \cdots
\sigma_{\gamma_{j-1}} (\gamma_j)$ for $ 1 \leq j \leq l$.
\end{prop}

\bigskip
 
\noindent
{\bf Proof.} For each $j = 1, 2, ..., l $, we have
\[
(t_j n_j t_{j}^{-1}) \circ \dot{\gamma}_j ~ = ~ t_j (n_j \circ
\dot{\gamma}_j) (\dot{\gamma}_j ~ t_{j}^{-1} ~  \dot{\gamma}_{j}^{-1}) ~ = ~
t_j (n_j \circ \dot{\gamma}_j) t_{j+1}^{-1},
\]
where $t_{\sl +1}$ is defined to be equal to $(\dg \dgm \cdots \dot{\gamma}_1)
~ t ~(\dg \dgm \cdots \dot{\gamma}_1)^{-1} = \dw^{-1} t \dw$.
Thus,
\[
\prod_{j=1}^{l} ~(t_j n_j t_{j}^{-1}) \circ \dot{\gamma}_j ~= ~
t \left( \prod_{j=1}^{l} (n_j \circ \dot{\gamma}_j) \right)
\dw^{-1} t \dw.
\]
On the other hand, let $n = F_w(n_1, ..., n_{\sl}) \in N_w$. We have
\[
(t ~n ~t^{-1}) \circ \dw ~ = ~ t (n \circ \dw) (\dw^{-1} t \dw).
\]
It now follows from the definition of $F_w$ that
\[
F_w(t_1 n_1 t_{1}^{-1}, ~ \cdots~
t_{\sl} n_{\sl} t_{\sl}^{-1}) ~ = ~ t ~n ~t^{-1}.
\]
Write $t = \exp(i H) \in T$ for $H \in \fa$ and $n_j = \exp(z_j
\check{E}_{\gamma_j}) \in N_{\gamma_j}$. Then (\ref{eq_t-cor})
follows from the following calculation:
\beqa
t_j ~ n_j ~ t_{j}^{-1} & = & \exp(z_j ~ Ad_{t_j} \check{E}_{\gamma_j})\\
& = & \exp(z_j ~ e^{i \sigma_{\gamma_1} \sigma_{\gamma_2} \cdots
\sigma_{\gamma_{j-1}}(\gamma_j) (H)} \check{E}_{\gamma_j}) \\
& = & \exp (t^{\alpha_j} z_j ~ \check{E}_{\gamma_j}).
\eeqa
\qed

\section{The Bruhat-Poisson structure}
\label{sec_bruhat}

Recall that a Poisson structure on 
a manifold $M$ (\cite{we:local}) is a bivector field $\pi
$ on $M$
such that the bracket operation on the algebra $C^{\infty}(M)$ of smooth
functions on $M$ defined by
\[
 \{\phi, \varphi\} ~ = ~ \pi(d \phi, d\varphi),  
\hspace{.3in} \phi, \varphi \in
C^{\infty}(M)
\]
satisfies Jacobi's identity. 
The condition on $\pi$ is that
\[
[\pi, ~ \pi ] ~ = ~ 0,
\]
where
$[ ~ , ~ ]$ denotes the Schouten bracket on the space
of multivector fields on $M$ (\cite{ko:crochet}). 

The bivector field $\pi$ can also be regarded as the bundle map
\[
 \tpi: T^{*}M \longrightarrow TM: ~ (\tpi(\alpha), ~ \beta) 
~ = ~ \pi(\alpha, ~ \beta).
\]
When $\pi$ is of maximal rank ($M$ is then necessarily even dimensional), the
bundle map $\tpi$ is invertible, and the 
$2$-form $\omega$ on $M$ defined by $\omega(x, y) = 
\pi (\tpi^{-1}(x), ~ \tpi^{-1}(y))$ is closed and non-degenerate and is thus
a symplectic $2$-form. In general,
the image of $\pi$ defines a (generally singular) involutive distribution 
on $M$. It has integrable submanifolds which inherit symplectic 
structures \cite{we:local}. They are called the symplectic leaves of $\pi$
in $P$.  Therefore, symplectic
manifolds are special cases of Poisson manifolds and every Poisson
manifold is a disjoint union of symplectic manifolds.

\bigskip
The Bruhat-Poisson structure on $K/T$ comes from a Poisson structure
$\pi$ on $K$ defined by
\[
\pi ~ = ~ \Lambda^r ~ - ~ \Lambda^l,
\]
where 
\[
\Lambda ~ = ~ {\frac{1}{2}} \sum_{\alpha > 0} X_{\alpha}
\wedge Y_{\alpha} ~ \in ~ \fk \wedge \fk
\]
with $\Xa$ and $\Ya$ given by (\ref{eq_xaya}), and $\Lambda^r$
(resp. $\Lambda^l$) is the right (resp. left) invariant
bi-vector field on $K$ with value $\Lambda$ at the identity element $e$.
We summerize some  properties of $\pi$ in the next theorem. For details,
see \cite{sts:dressing} \cite{soi:compact} and
\cite{lu-we:poi}.

\bigskip
\begin{thm}
\label{thm_on-k}

a) The bi-vector field $\pi$ defines a Poisson structure on $K$;

b) Equip $K \times K$ with the product Poisson structure $\pi \oplus \pi$. 
Then
the multiplication map 
\[
K \times K \lrw K: ~~ (k_1, ~ k_2) \Map k_1 k_2
\]
is a Poisson map, making $(K, \pi)$ into a Poisson Lie group;

c) The symplectic leaves of $\pi$ in $K$ are precisely the 
orbits of $AN$ in $K$ for the action given by (\ref{eq_G-on-K}).
These are the Bruhat cells in $K$: for each Weyl group
element $w \in W$ with a fixed representative $\dw$ of $w$ in $K$, 
the symplectic leaf through $\dw$ is the $N_w$-orbit 
$C_{\dot{w}} = N_{w} \circ \dw$ introduced in Section \ref{sec_cor}.
For each $t \in T$, the subset $C_{\dot{w}} t$ is the symplectic leaf
through the point $\dw t$. As $w$ runs over $W$ and $t$ over $T$, these 
are all the symplectic leaves in $K$;

d) Both left and right translations by elements in $T$ leave $\pi$ invariant;

e) The image of $\pi$ under the projection
map $K \rightarrow K/T$ is a well defined bi-vector field on $K/T$
which we still denote by $\pi$. It defines a Poisson structure on
$K/T$ called the Bruhat-Poisson structure;

f) Symplectic leaves of the Bruhat-Poisson structure on $K/T$
are precisely the Bruhat cells $\Sigma_w$, for
$w \in W$, in $K/T$ (and thus the name);

g) With respect to the left translations by elements in $K$, the 
Bruhat-Poisson structure is $T$-invariant but not $K$-invariant;
The action map
\[
K \times K/T \lrw K/T: ~~ (k_1, ~ k_2/T) \Map k_1 k_2/T
\]
is a Poisson map, making $(K/T, ~ \pi)$ into a Poisson
homogeneous $(K, ~ \pi)$-space.
\end{thm}

\bigskip
Therefore, for each $w \in W$ with a fixed representative 
$\dw$ in $K$, both $C_{\dot{w}}$ and $\Sigma_w$ inherit
symplectic structures 
as symplectic leaves in $K$ and $K/T$
respectively, and the projection from $C_{\dot{w}}$ to 
$\Sigma_w$ is a symplectic diffeomorphism. Moreover, the symplectic structure
on $\Sigma_w$ is invariant under the action of $T$ by left translations. 
The goal of this section is to 
write down both the symplectic structure - 
thus also the Liouville measure - on $\Sigma_w$ and the moment
map for the $T$-action in the 
$\zc$ coordinates. Here, we regard 
$\zc$ as coordinates on $\Sigma_w$ via the parametrization
of $\Sigma_w$ by $N_w$,

Property b) of $\pi$ on $K$ is called the multiplicativity of $\pi$. 
As we will see shortly, it is exactly this property that enables us
to decompose, as symplectic manifolds, the higher dimensional 
$C_{\dot{w}}$'s into products of $2$-dimensional ones.
In fact, this is also the motivation for introducing the coordinates
$\zc$.

\bigskip
As in Section \ref{sec_cor}, let 
$l = l(w)$ and  fix a  reduced decomposition
\[
w ~ = ~ \sigma_{\gamma_1} \sigma_{\gamma_2} \cdots 
\sigma_{\gamma_l}.
\]
For each $j = 1, ..., l$, let
$\dot{\gamma}_j$ be a representative of $\sigma_{\gamma_j}$
in $K$, so $\dw = \dc_1 \dc_2 \cdots \dc_{\sl}$ is a
representative of $w$ in $K$. Property c) of $\pi$ on $K$
says that the symplectic leaf through $\dc_j$ is 
the $2$-dimensional cell $C_{\dot{\gamma}_j}$. Theorem 
\ref{thm_cor} (see Remark \ref{rem_soibel}) says that the map
\begin{equation}
\label{eq_c}
C_{\dot{\gamma}_1}  \times C_{\dot{\gamma}_1}  \times \cdots 
\times C_{\dot{\gamma}_{\sl}} \lrw C_{\dot{w}}.
\end{equation}
is a diffeomorphism.

\begin{prop}
\label{prop_c}
The map in (\ref{eq_c}) is a symplectic diffeomorphism, where the left
hand side has the product symplectic structure.
\end{prop}

\noindent
{\bf Proof.} 
(Notice how the multiplicativity of $\pi$ is used in the proof.)
The inclusion map 
\[
C_{\dot{\gamma}_1}  \times C_{\dot{\gamma}_1}  \times \cdots 
\times C_{\dot{\gamma}_{\sl}} ~ \hookrightarrow ~
K \times K \times \cdots \times K \hspace{.3in}
(\text{ $l$-copies})
\]
is a Poisson map because each $C_{\dot{\gamma}_j}$ is a symplectic leaf
and thus a Poisson submanifold of $K$. The multiplicativity of $\pi$
says that the multiplication map
\[
K \times K \times \cdots \times K  \lrw K: ~~
(k_1, ~ k_2, ~ \cdots, ~ k_{\sl}) \Map k_1 k_2 \cdots k_{\sl}
\]
is a Poisson map. Thus, composing the two, we see that 
\[
C_{\dot{\gamma}_1}  \times C_{\dot{\gamma}_1}  \times \cdots 
\times C_{\dot{\gamma}_{\sl}} \lrw 
K: ~~ (k_1, ~ k_2, ~ \cdots, ~ k_{\sl}) \Map k_1 k_2 \cdots k_{\sl}
\]
is a Poisson map. But it has its image in $C_{\dot{w}}$ which is
 a Poisson submanifold of $K$. Thus, regarded as a map to 
$C_{\dot{w}}$, the above is a Poisson map, and thus a Poisson 
and therefore a symplectic diffeomorphism.   
\qed

It now remains to determine the symplectic structure on the 
$2$-dimensional leaves.  Recall that for each simple root
$\gamma$, 
 we have a Lie group homomorphism
\[
\Phi_{\gamma}: ~~ SL(2, ~ {\Bbb C}) \lrw G
\]
which maps $SU(2)$ to $K$. 

\begin{prop}
\label{prop_su2}
(See also \cite{soi:compact}) 
For each simple root $\gamma$, equip $SU(2)$ with the Poisson
structure given by
\[
\pi_{\gamma} ~ = ~ \Lambda_{\gamma}^{r} ~ - ~ \Lambda_{\gamma}^{l},
\]
where 
\[
\Lambda_{\gamma} ~ = ~ {\frac{1}{4}} \ll \gamma, ~ \gamma \gg 
\left( \begin{array}{cc} 0 & 1 \\ -1 & 0 \end{array} \right) 
\wedge 
\left( \begin{array}{cc} 0 & i \\ i & 0 \end{array} \right) 
~ \in ~ su(2) \wedge su(2),
\]
and $\Lambda_{\gamma}^{r}$ (resp. $\Lambda_{\gamma}^{l}$)
denotes the right (resp. left) invariant bi-vector field 
on $SU(2)$ with value $\Lambda_{\gamma}$ at the identity element.
Then 1) $\Psi_{\gamma}: ~ (SU(2), ~ 
\pi_{\gamma}) \rightarrow
(K, ~ \pi)$
is a Poisson map, and 2)
the symplectic leaf of $\pi_{\gamma}$ in $SU(2)$ through the 
point $\left( \begin{array}{cc} 0 & i \vspace{-.05in}
\\ i & 0 \end{array} \right)$
is 
\[
C_0 ~ = ~ \left\{ 
\left( \begin{array}{cc}
{\dfrac{i z}{\sqrt{1 + |z|^2}}} & {\dfrac{i}{\sqrt{1 + |z|^2}}}
\vspace{.1in}\\
{\dfrac{i}{\sqrt{1 + |z|^2}}} & {\dfrac{-i \bar{z}}{\sqrt{1 + |z|^2}}}
\end{array} \right): ~~~ z  \in {\Bbb C}
\right\}.
\]
Using $(z, ~ \bar{z})$ as coordinates on $C_0$, the induced 
symplectic structure is given by
\[
\Omega ~ = ~ {\frac{i}{\ll \gamma, ~ \gamma \gg}}~ 
{\frac{1}{1+|z|^2}} dz \wedge d \bar{z}.
\]
\end{prop}

\bigskip
\noindent
{\bf Proof.} 1) Think of $\fa + \fn$ as a real Lie algebra and identify 
$\fa + \fn$ with $\fk^*$ using the
imaginary part of the Killing form as the pairing. 
The fact that $\Psi_{\gamma}$ is a Poisson map then follows from the fact
that the subspace $\Psi_{\gamma}(su(2))^{\perp}$ of $\fa + \fn$
consisting of all elements that annihilate 
$\Psi_{\gamma}(su(2))$ is an ideal of
$\fa + \fn$. See \cite{lu-we:poi}. 

2) As a special case of Theorem \ref{thm_on-k} (up to a constant
multiple), 
the symplectic leaves of $\pi_{\gamma}$ in $SU(2)$ are either $2$ or
$0$-dimensional. We know from Section \ref{sec_cor} that the set $C_0$ 
is the symplectic leaf through the point 
$\left( \begin{array}{cc} 0 & i \vspace{-.05in}\\ i & 0 \end{array} \right)$.

Write an element of $SU(2)$ as 
\[
u ~ = ~ \left( \begin{array}{cc} \xi & \eta \vspace{.05in} \\ -\bar{\eta} & 
\bar{\xi} \end{array} \right).
\]
The Poisson brackets defined by $\pi_{\gamma}$ are
\beqa
\{ \xi, ~ \bar{\xi} \} & = & -i \ll \gamma, ~ \gamma \gg |\eta|^2 
\hspace{.4in}
 \{ \xi, ~ \eta \} ~= ~ {\frac{1}{2}} i \ll \gamma, ~ \gamma \gg
\xi \eta \\
\{ \xi, ~ \bar{\eta} \} & = & {\frac{1}{2}} i  \ll \gamma, ~ \gamma \gg
\xi \bar{\eta} \hspace{.5in} 
\{ \eta, ~ \bar{\eta} \} ~ = ~ 0.
\eeqa
Set
\[
\xi ~ = ~ {\frac{iz}{\sqrt{1+|z|^2}}}, \hspace{.4in}
\eta ~ = ~ {\frac{i}{\sqrt{1+|z|^2}}},
\]
so
\[
z ~ = ~ \xi / \eta, \hspace{.4in} 
\bar{z} ~ = ~ \bar{\xi} / \bar{\eta}.
\]
Using the Leibniz rule for the Poisson bracket, we get
\beqa
\{ z, ~ \bar{z} \} & = & \{\xi / \eta, ~~ \bar{\xi} / \bar{\eta} \} \\
& = & {\frac{1}{|\eta|^2}} \{\xi, ~ \bar{\xi} \} 
~ - ~ {\frac{\bar{\xi}}{|\eta|^2 \bar{\eta}}} 
\{ \xi, ~ \bar{\eta} \} ~ - ~ 
{\frac{\xi}{|\eta|^2 \eta}} \{\eta, ~ \bar{\xi} \} \\
& = & -i \ll \gamma, ~ \gamma \gg (1 ~ + ~ {\frac{|\xi|^2}{|\eta|^2}})\\
& = & -i \ll \gamma, ~ \gamma \gg (1 + |z|^2).
\eeqa
Thus the symplectic $2$-form $\Omega$ on $C_0$ is given as 
stated.
\qed

\begin{thm}
\label{thm_sym-cor}
In the $\zc$ coordinates, the induced symplectic structure on 
$\Sigma_w$ (as a symplectic leaf of the Bruhat-Poisson structure) 
and the Liouville volume form $\mu_w$ associated to $\Omega_w$
are given by
\begin{eqnarray}
\label{eq_Omega}
\Omega_w & = & \sum_{j = 1}^{l} 
{\frac{i}{\ll \alpha_j, ~ \alpha_j \gg}} 
~ {\frac{1}{1 + |z_j|^2}} ~dz_j \wedge d \bar{z}_j\\
\label{eq_mu-cor}
\mu_w &: = &  {\frac{1}{l!}} ~ (\Omega_{w})^{l}
~ = ~ \prod_{j=1}^{l} {\frac{i}{\ll \alpha_j, ~ \alpha_j \gg}}
~ {\frac{1}{1 + |z_j|^2}} ~ dz_j \wedge d \bar{z}_j.
\end{eqnarray}
The moment map for the $T$-action on $\Sigma_w$ by left 
translations satisfying $\phi_w(0) = 0$ is given by
\begin{equation}
\label{eq_mom}
\phi_w: ~~ \Sigma_w \lrw \ft^* \cong \fa: ~~ \phi_w ~ = ~ 
\sum_{j=1}^l (-{\frac{1}{2}} \log(1+|z_j|^2) \check{H}_{\alpha_j}).
\end{equation}
Here we are identifying $\ft^*$ with $\fa$ by using the imaginary part of
the Killing form as the pairing.
\end{thm}

\noindent
{\bf Proof.} The formula for $\Omega_w$ follows immediately from 
Propositions \ref{prop_c} and \ref{prop_su2} and the fact that
$\ll \gamma_j, ~ \gamma_j \gg = \ll \alpha_j, ~ \alpha_j \gg$.

The formula for the moment $\phi_w$ follows 
from the explicit formula for the $T$-action on $\Sigma_w$
in the $\zc$ coordinates as given in Proposition \ref{prop_t-on-N}: let
$iH \in \ft$, where $H \in \fa$. By Proposition \ref{prop_t-on-N},
the generating vector field of the $T$ action on $\Sigma_w$ 
in the direction of $iH$ is given by
\[
V_{iH} ~ = ~ \sum_{j = 1}^{l} \alpha_j(H) (-y_j 
{\frac{\partial}{\partial x_j}} + x_j 
{\frac{\partial}{\partial y_j}}).
\]
Set
\[
V_{iH} \backl \Omega_w ~ = ~ d < \phi_w, ~ iH >.
\]
We see that the only solution of 
$\phi_w$ that satisfies this equation for all
$H \in \fa$ and the condition that $\phi_w(0) = 0$ is
given by (\ref{eq_mom}).
\qed

The following two corollaries follow immediately by comparing
the formulas in Theorems
\ref{thm_a}, \ref{thm_haar} and \ref{thm_sym-cor} (see also identity 
(\ref{eq_a-w})).

\begin{cor}
\label{cor_aw}
Think of the moment map $\phi_w$ as a map from $N_w$ to $\fa$ via
the parametrization of $\Sigma_w$ by $N_w$. Then 
\[
\phi_w ~ = ~ Ad_{\dot{w}} \log a_w(n).
\]
\end{cor}

\begin{cor}
\label{cor_liouville-haar}
Think of the Haar measure $dn$ of $N_w$ as a volume form on $\Sigma_w$
via the parametrization of $\Sigma_w$ by $N_w$. It is related to 
the Liouville measure $\mu_w$ by
\[
\mu_w ~ = ~ \left(\prod_{j=1}^{l} {\frac{\pi}{\ll \rho, ~ \beta_j \gg}} 
\right) ~ a_{w}(n)^{-2 \rho} dn.
\]
\end{cor}

\bigskip
We have thus connected the moment map $\phi_w$ and the Liouville measure
$\mu_w$ with the familiar map 
$a_w: N_w \rightarrow A$ and the Haar measure $dn$ on $N_w$. 
Such a connection is desirable for understanding the geometry of the 
Bruhat-Poisson structure.
We have arrived at this by comparing their formulas
in coordinates. This is why we wanted to write down the
formulas for $a_w$ and $dn$ in the $\zc$ coordinates in Section \ref{sec_cor}.
 
\section{Kostant's theorem on $H(G/B)$}
\label{sec_kost-thm}

In \cite{ks:63}, Kostant
constructs, for each element $w$ in the Weyl group $W$, an explicit
$K$-invariant closed differential form $s^{w}$ on $X$ with
$\deg(s^w) = 2 l(w)$, such that the cohomology classes of the $s^w$'s 
form a basis of $H(X, {\Bbb C})$
that, up to scalar multiples, is dual to the basis of the 
homology of $X$ formed by the
closures of the Bruhat cells in $X$. We now recall the definition 
of these forms in more details. 

\bigskip
Let $C = \oplus C^{\bullet}$ be the space of all $K$-invariant 
differential forms on $G/B \cong K/T$. By identifying $\fg$ and $\fg^*$
via the Killing form $\ll ~ , ~ \gg$ so that
\[
(\fg /\fh)^* ~ \cong ~ \fn_{-} ~ + ~ \fn,
\]
we can identify 
\[
C^{\bullet} ~ \cong ~ \wedge^{\bullet} (\fn_{-} ~ + ~ \fn)^{T}.
\]
Introduce the operators $E$ and $L_0$ on $C^{\bullet}$ by
\[
E ~ = ~ 2 \sum_{\alpha > 0} ad_{E_{-\alpha}} ~ \ot ~ 
ad_{E_{\alpha}}
\]
and
\beqa
& & L_0(E_{-\alpha_1} \wedge E_{\alpha_2} \wedge \cdots \wedge E_{-\alpha_{p}}
\ot  E_{\beta_1} \wedge E_{\beta_2}  \wedge \cdots \wedge 
E_{\beta_{q}}) \vspace{.2in} \\
 \hspace{.4in}  &  = & \left\{
\begin{array}{ll} 0 & \text{if} ~ \| \rho \|^2 - \|\rho -(\alpha_1 
+ \alpha_2 + \cdots + \alpha_p) \|^2 = 0 \vspace{.1in} \\
{\dfrac{1}{ \| \rho \|^2 - \|\rho -(\alpha_1  
+ \alpha_2 + \cdots + \alpha_p) \|^2}} & \text{if} ~ 
 \| \rho \|^2 - \|\rho -(\alpha_1  
+ \alpha_2 +  \cdots + \alpha_p) \|^2 \neq 0 \end{array} \right.
\eeqa
Set \[
R ~ = ~ - L_0 E.
\]
It is a nilpotent operator. 

\bigskip
Now let $w \in W$ be a Weyl group element. As before, use $l$
to denote $l(w)$. Let 
\[
\{ \alpha_1, ~ \alpha_2, ~ ..., ~ \alpha_{\sl} \} ~ = ~ 
\{ \alpha > 0: ~~ w^{-1} \alpha < 0 \},
\]
and let
\[
\beta_j ~ = ~ - w^{-1} \alpha_j
\]
so that
\[
\{ \beta_1, ~ \beta_2, ..., ~ \beta_{\sl} \} ~ = ~ \{ \beta > 0: ~~
w \beta < 0 \}.
\]
Set
\[
h^{w^{-1}} ~ = ~ ({\frac{i}{2}})^{l} ~ E_{-\beta_1} \wedge E_{-\beta_2}
\wedge \cdots \wedge E_{-\beta_{\sl}} ~ \ot ~ 
E_{\beta_1} \wedge E_{\beta_2} \wedge \cdots \wedge E_{\beta_{\sl}},
\]
and
\begin{equation}
\label{eq_s}
s^w ~ = ~ (1 - R)^{-1} h^{w^{-1}} ~ = ~ h^{w^{-1}} ~ + ~ 
R h^{w^{-1}} ~ + ~ R^2 h^{w^{-1}} ~ + ~ \cdots.
\end{equation}

\bigskip
\begin{thm} [Kostant \cite{ks:63}]
\label{thm_kost}
1). The forms $s^{w}$ for $w \in W$ are closed, and their cohomology classes
form a basis of the de Rham cohomology of $G/B \cong K/T$ (with 
complex coefficients);

2). Let $j_w: ~ N_w \rightarrow \Sigma_w: n \mapsto n w/T$ be
the parametrization map. Then 
\[
j_{w_1}^{*} (s^w|_{\Sigma_{w_1}}) ~ = ~ 
\left\{ \begin{array}{ll} 0 & \text{if} ~ l(w_1) = l(w) ~ 
\text{but} ~w_1 \neq w \\ a_w(n)^{-2(\rho - w^{-1} \rho)} (dn)_1
& \text{if} ~ w_1 = w, \end{array} \right.
\]
where $s^w|_{\Sigma_{w_1}} = i_{w_1}^{*} s^w$ if $i_{w_1}$ is the 
inclusion map of $\Sigma_{w_1}$ into $K/T$, the map $a_w: 
N_w \rightarrow A$ is, 
as before, defined by (\ref{eq_aw}), and
$(dn)_1$ is the left invariant Haar measure on $N_w$ normalized
by the condition
\[
\left( (dn)_1, ~ E_{\alpha_1} \wedge i E_{\alpha_1} 
\wedge E_{\alpha_2} \wedge i E_{\alpha_2} \wedge \cdots \wedge 
E_{\alpha_{\sl}} \wedge i E_{\alpha_{\sl}} \right) ~ = ~ 1.
\]

3)
\[
\int_{\Sigma_w} s^w ~ = ~ \prod_{j=1}^{l}
{\frac{\pi}{ \ll \rho, ~ \alpha_j \gg}}.
\]
\end{thm}

\bigskip
The purpose of this section is to relate the form $s^w$ with
the Liouville volume form $\mu_w$ on $\Sigma_w$ induced by
the Bruhat-Poisson structure. 
This is now easy due to Corollaries
\ref{cor_aw} and \ref{cor_liouville-haar}.
We will also give a simple proof of 3).

\bigskip
We first need to relate the left Haar measures $(dn)_1$ and
$dn$ as given in Theorem \ref{thm_haar}. 

\begin{lem}
\label{lem_haar}
We have
\[
(dn)_1 ~ = ~ \prod_{j=1}^{l}
{\frac{\pi}{ \ll \rho, ~ \beta_j \gg}} dn.
\]
\end{lem}

\noindent
{\bf Proof.} This follows from Proposition \ref{prop_scalar} and
the fact that $\ll \alpha_j, ~ \alpha_j \gg = \ll \beta_j, ~ \beta_j \gg$
for each $j = 1, ..., l$.
\qed
  
\begin{thm}
\label{thm_main}
When restricted to the Schubert cell $\Sigma_w$, Kostant's harmonic 
form $s^w$ is related to the Liouville volume form $\mu_w$ on $\Sigma$ by
\begin{equation}
\label{eq_main}
s^w|_{\Sigma_w} ~ = ~ (a_{w})^{2 w^{\scriptscriptstyle -1} \rho} \mu_w ~ = ~ 
 e^{<\phi_w, ~ 2 i H_{\rho} >} \mu_w.
\end{equation}
\end{thm}

\bigskip
\noindent
{\bf Proof.} 
This is a direct consequence of Corollaries \ref{cor_aw} and 
\ref{cor_liouville-haar}. Explicitly, 
we have
\begin{eqnarray}
\label{eq_s-cor}
& & s^w|_{\Sigma_w} ~ = ~ \prod_{j=1}^{l} 
{\frac{i }{\ll \alpha_j, ~ \alpha_j \gg}} (1+|z_j|^2)^{
{-~ \dfrac{2 \ll \rho, ~ \alpha_j \gg}{\ll \alpha_j, ~ \alpha_j \gg}} 
~ -~ 1}
dz_j \wedge d \bar{z}_j\\
\label{eq_mom-cor}
& & <\phi_w, ~ 2 i H_{\rho} > ~ = ~ - \sum_{j=1}^{l} 
{\frac{2\ll \rho, ~ \alpha_j \gg}{\ll \alpha_j, ~ \alpha_j \gg}}
\log(1 + |z_j|^2)
\end{eqnarray}
and
\[
\mu_w ~ = ~ \prod_{j=1}^{l} {\frac{i}{\ll \alpha_j, ~ \alpha_j \gg}}
(1+|z_j|^2)^{-1} dz_j \wedge d \bar{z}_j.
\]
Thus we have (\ref{eq_main}).
\qed

From the explicit formula for $s^w|_{\Sigma_w}$,
we immediately get
\beqa
\int_{\Sigma_w} s^w  & = &  \prod_{j=1}^{l} \int_{{\Bbb R}^2}
{\frac{i}{\ll \alpha_j, ~ \alpha_j \gg}} (1+|z_j|^2)^{
{-~ \dfrac{2 \ll \rho, ~ \alpha_j \gg}{\ll \alpha_j, ~ \alpha_j \gg}} 
~ - ~1}
dz_j \wedge d \bar{z}_j \\
& = & \prod_{j=1}^{l} {\frac{\pi}{ \ll \rho, ~ \alpha_j \gg}}.
\eeqa
This integral was first calculated in \cite{kk:integral} using
induction on $l(w)$. Again, as in the case of the $c$-functions,
our simple proof is due to the fact that the induction
argument has been pushed to the calculations for $a_w(n)$ and $dn$.
This is in fact a special case of the $c$-function with
$i \lambda = -2 w^{-1} \rho$.

\begin{rem}
\label{rem_haar-on-k}
{\em
When $w = w_0$ is the longest Weyl group element, the form $s^{w_0}$
is a $K$-invariant volume form so it coincides with the Haar measure
on $K/T$. When restricted to the biggest cell $\Sigma_{w_0}$, the
Liouville volume form $\mu_{w_0}$, the form $s^{w_0}$, and the
Haar measure $dn$ for $N_{w_0} = N$ are related by 
\[
s^{w_0}|_{\Sigma_{w_0}} ~ = ~ (a_{w_0})^{-2 \rho} \mu_{w_0} 
~ = ~ \left(\prod_{\alpha > 0} {\dfrac{\pi}{\ll \rho, ~ \alpha \gg}} \right)
 a_{w_0}(n)^{-4 \rho} dn.
\]
}
\end{rem}

\begin{rem}
\label{rem_not-hamil}
{\em
The function $<\phi_w, ~ 2iH_{\rho} >$ on $\Sigma_w$ is
the Hamiltonian function for the
generating vector field  $\theta_0$ of the $T$-action in the direction of
$2 i H_{\rho}$. This vector field is intrinsic to the Bruhat-Poisson
structure in the sense that it is its modular vector field. 
As $|z_j| \rightarrow +\infty$, the function $<\phi_w, ~ 2 i H_{\rho} > 
\rightarrow - \infty$. Thus the modular vector field $\theta_0$ is not
globally Hamiltonian on $K/T$. We say that the Bruhat-Poisson structure
is not unimodular. See \cite{we:modular} \cite{bz:outer} \cite{elw:modular}.
}
\end{rem}

\begin{rem}
\label{rem_sam}
{\em
The second identity in (\ref{eq_main}) expresses the form $s^w|_{\Sigma_w}$
totally in terms of data coming from the Bruhat-Poisson structure.
In particular, it says that the integral $\int_{\Sigma_w} s^w$
is of the Duistermaat-Heckman type. Wanting to see this was another
motivation for this work.  Theorem \ref{thm_main} can be used to 
describe generators of the so-called Poisson - de Rham cohomology
of the Bruhat-Poisson structure. We do this in
\cite{e-l:poi}. 
}
\end{rem}

\section{Appendix: Relation to the Bott-Samelson coordinates}
\label{sec_bs}

As before, let $w \in W$ be a Weyl group element with a fixed
reduced decomposition:
\[
w ~ = ~ \sigma_{\gamma_1} \sigma_{\gamma_2} \cdots
\sigma_{\gamma_{l}}
\]
where $l = l(w)$. Then we have the Lie group homomorphism
$\Psi_{\gamma_j}: SL(2, {\Bbb C}) \rightarrow G$ and the element
\[
\dot{\gamma}_j ~ = ~ \Psi_{\gamma_j}
\left( \begin{array}{cc} 0 & i \\ i & 0 \end{array} \right) ~ \in~ K
\]
for each $j = 1, ..., l$. Set again
$\dot{w} = \dot{\gamma_1} \dot{\gamma_2} \cdots \dot{\gamma_{\sl}} \in K$.
Then the map
\beqa
& & F_{w}': ~ N_{\gamma_1} \times \cdots \times N_{\gamma_{\sl}} 
  \lrw  N_w: \\ 
& & (n_1, ~ n_2, ~ ... ~ , n_{\sl})  \Map  n_1 \dot{\gamma}_1 
n_2 \dot{\gamma}_2 \cdots n_{\sl}\dot{\gamma}_{\sl} \dw^{-1} \\
& & \hspace{1.3in} = ~  n_1 ~ (\dot{\gamma}_1 n_2 \dot{\gamma}_{1}^{-1}) 
\cdots (\dot{\gamma}_1 \dot{\gamma}_2 \cdots \dgm) ~ n_{\sl} ~
(\dot{\gamma}_1 \dot{\gamma}_2 \cdots \dgm)^{-1}
\eeqa
is a holomorphic diffeomorphism.
Thus, if we parametrize $N_{\gamma_j}$ by ${\Bbb C}$ via
\[
 z_j' \Map n({z_j'}) ~ := ~ \Psi_{\gamma_j} \left( 
\begin{array}{cc} 1 & z_j' \\ 0 & 1 \end{array} \right)
~ = ~ \exp(z_j' \check{E}_{\gamma_j}) ~ \in ~ N_{\gamma_j}
\]
for each $j$, we get a parametrization of $N_w$
by ${\Bbb C}^l$: 
\[
{\Bbb C}^{l} \lrw N_w: ~~ (z_1', ~ z_2', ~ ... ~ , z_{\sl}')
\Map F_w' (n({z_1'}), ~ n({z_2'}), ~ ... ~ , n({z_{\sl}'})).
\]
We call it the Bott-Samelson
parametrization of $N_w$ and call $\{z_1', ~ z_2', ~ ... ~ , z_{\sl}'\}$ 
the {\bf Bott-Samelson coordinates} on $N_w$ because of the close relation
to
the Bott-Samelson desingularization of the Schubert variety $X_w$ 
\cite{j:group}.

\bigskip 
The difference between our coordinates $\zc$ and the Bott-Samelson 
coordinates $\{z_1', ~ z_2', ~ ... ~ , z_{\sl}'\}$ can
be described as follows: by composing  $F_w$ (in Theorem \ref{thm_cor})
and $F_w'$ with the map
\[
j_w: ~~ N_w \lrw \Sigma_w: ~~ n \Map n \dw / B,
\]
we can think of both $F_w$ and $F_w'$ as mapping 
$N_{\gamma_1} \times \cdots \times N_{\gamma_{\sl}}$ diffeomorphically 
to $\Sigma_w$. 
The map $F_w'$ first sends $(n_1, n_2, ..., n_{\sl})$ to the product
of $n_1 \dot{\gamma}_1  n_2 \dot{\gamma}_2 \cdots  n_{\sl}
\dot{\gamma}_{\sl}$ in $G$ and then projects the product to 
$G/B$. But for 
$F_w$, we first pick up the $K$-component
$k_j$ in the Iwasawa decomposition of $n_j \dot{\gamma}_j$ for
each $j$, 
multiply the $k_j$'s inside $K$ and then project the product to $K/T 
\cong G/B$. 

Consider now the change of coordinates
\[
(z_1', ~ z_2', ~ ..., ~ z_{\sl}') ~ = ~ I_w (z_1, ~ z_2, ..., ~ z_{\sl}).
\]
We can get a recursive formula for $I_w$ from the proof of Theorem 
\ref{thm_haar}. Clearly, when $w$ is a simple reflection, the
map $I_w$ is the identity map. For a general $w$ with
the reduced decomposition  $w= 
\sigma_{\gamma_1} \sigma_{\gamma_2} \cdots \sigma_{\gamma_{\sl}}$, 
let again $w_1 =
\sigma_{\gamma_1} \sigma_{\gamma_2} \cdots \sigma_{\gamma_{\sl-1}}$
so that $w = w_1 \sigma_{\gamma_{\sl}}$. Keeping the same notation as
in the proof of Theorem \ref{thm_haar}, we know that $I_w$ is
given by
\[
\left\{ \begin{array}{l} 
(z_1', ~ z_2', ~ ..., ~ z_{\sl -1}') ~ = ~ I_{w_1} 
(z_1, ~ z_2, ..., ~ z_{\sl-1}) \vspace{.1in} \\
z_{\sl}' ~ = ~ m(z_1, ~ z_2, ..., ~ z_{\sl-1}) ~ + ~ 
a_{w_1}(n)^{-\gamma_{\sl}} z_{\sl}, \end{array} \right.
\]
where, recall, $ m(z_1, ~ z_2, ..., ~ z_{\sl-1}) \in {\Bbb C}$
is such that $\exp( m(z_1,  z_2, ...,  z_{\sl-1}) 
\check{E}_{\gamma_{\sl}})$ is the $N_{\gamma_{\sl}}$-component of
$(m')^{-1} \in N$ with respect to the decomposition
$N = N_{\gamma_{\sl}} \Nhat$, and $m' \in N$ is the $N$-component 
in the Iwasawa decomposition of $n \dw_1$. 
The change of coordinates is in general not complex because
the function $m(z_1, ~ z_2, ..., ~ z_{\sl-1})$ is in general 
not holomorphic. This is also seen from the 
following example.

\begin{exam}
\label{exam_sl3-again}
{\em
For $G = SL(3, {\Bbb C})$ and $w = (1,2)(2,3)(1,2)$, the Bott-Samelson
parametrization of $N_w$ is by
\[
(z_1', ~ z_2', ~ z_3') \Map \left( \begin{array}{ccc} 
1 & z_1' & z_1' z_3' + i z_2' \hspace{.1in}\\ 0 & 1 & z_3' 
\hspace{.1in}\\ 0 & 0 & 1 \end{array}
\right).
\]
The change of
coordinates between $(z_1, ~ z_2, ~ z_3)$
and $\{z_1', ~ z_2', ~ z_3'\}$
are (see Example \ref{exam_sl3})
\[
z_1' ~ = ~ z_1, \hspace{.4in} 
z_2' ~ = ~ \epsilon_1 z_2 \hspace{.4in}
z_3' ~ = ~ {\dfrac{\epsilon_2 z_3 - i \bar{z}_1 z_2}{\epsilon_1}},
\]
or
\[
z_1 ~ = ~ z_1', \hspace{.4in}
z_2 ~ = ~ {\dfrac{z_2'}{\eta_1}} \hspace{.4in}
z_3 ~ = ~ {\dfrac{\eta_{1}^{2} z_3' + i \overline{z_1'} z_2'}{\eta_2}},
\]
where $\epsilon_j = \sqrt{1 + |z_j|^2}$ for $j = 1, 2$ and 
$ \eta_1  =  \sqrt{1  +  |z_1'|^2},  ~ 
\eta_2  =  \sqrt{1  +  |z_1'|^2  + |z_2'|^2}$.
}
\end{exam}

\end{document}